\documentclass[11pt]{article}%
\usepackage{algpseudocode}
\usepackage{float} % fix table inside context 
\usepackage{fullpage}
\usepackage{amsfonts}
\usepackage{amsmath}
\usepackage{amssymb}
\usepackage{amsbsy}
\usepackage{amsthm}
\usepackage{mathtools}
\usepackage{epsfig}
\usepackage{graphicx}
\usepackage{colordvi}
\usepackage{graphics}
\usepackage{color}
\usepackage{bm}
\usepackage{verbatim} 
\usepackage{blkarray}
\numberwithin{equation}{section}

\newtheorem{theorem}{Theorem}[section]
\newtheorem{thm}{Theorem}[section]
\newtheorem{algorithm}{Algorithm}[section]
\newtheorem{cor}[thm]{Corollary}
\newtheorem{lemma}[thm]{Lemma}

\newtheorem{prop}[thm]{Proposition}
\newtheorem{remark}[thm]{Remark}
\newtheorem{assumption}[thm]{Assumption}

\newcommand{\R}{\mathbb{R}}
\newcommand{\C}{\mathbb{C}}
\newcommand{\norm}[1]{\left\Vert#1\right\Vert}
\newcommand{\nr}[1]{\ensuremath{\left\|{#1}\right\|}}

\newcommand{\anorm}{\norm{\,\cdot \,}}

\DeclareMathOperator{\diag}{diag}
\newcommand{\eps}{\varepsilon}

\newcommand{\goto}{\rightarrow}

\newcommand{\cP}{{\mathcal P}}
\newcommand{\bigo}{{\mathcal O}}

\newcommand{\f}{\frac}

\DeclareMathOperator{\sign}{sign}
\DeclareMathOperator{\spa}{span}
\DeclareMathOperator{\dist}{dist}

\begin{document}

\title{Dynamically accelerating the power iteration with momentum}

\author{
Christian Austin
\thanks{Department of Mathematics, University of Florida, Gainesville, FL
  32611-8105 (christianaustin@ufl.edu)}
\and 
Sara Pollock
\thanks{Department of Mathematics, University of Florida, Gainesville, FL
  32611-8105 (s.pollock@ufl.edu)}
\and
Yunrong Zhu
\thanks{Department of Mathematics \& Statistics, Idaho State University,
  Pocatello, ID 83209 (zhuyunr@isu.edu)}
}
\date{}
%\subkclass[2010]{}
%\keywords{
%}

\maketitle

\begin{abstract}
In this paper, we propose, analyze and demonstrate 
a dynamic momentum method to accelerate power and 
inverse power iterations with minimal computational overhead. 
The method can be applied to real diagonalizable matrices, 
is provably convergent with acceleration in the symmetric case,
and does not require a priori spectral knowledge.
We review and extend background results on previously developed static momentum 
accelerations for the power iteration through the connection between the momentum 
accelerated iteration and the standard power iteration applied to an augmented matrix.
We show that the augmented matrix is defective for the optimal parameter choice.
We then present our dynamic method which updates the momentum parameter at each
iteration based on the Rayleigh quotient and two previous residuals. 
We present convergence and stability theory for the method by considering a power-like
method consisting of multiplying an initial vector by a sequence of augmented matrices.
We demonstrate the developed method on a number of benchmark problems, and see that
it outperforms both the power iteration and often the static momentum 
acceleration with optimal parameter choice. Finally, we present and
demonstrate an explicit extension of the algorithm to inverse power iterations.
\end{abstract}

%% -- --------------------------------------------------------------------
\section{Introduction}\label{sec:intro}

In recent years, there is a resurgence of interest in the power method, given its 
simplicity and ease of implementation. This method to find the dominant eigenmode 
of a matrix can be applied in a variety of machine 
learning algorithms, such as PCA, clustering, and low-rank matrix approximations 
(see \cite{RJRH22} and the references cited therein), PageRank \cite{BRZ06,GoGr06,HKMG03,kamvar10,sidi04}, 
and stability analysis of partial differential equations \cite{BaOs91}. 

There are a number of generalizations of the power method for large and often
sparse systems that 
can be used to compute extreme eigenvectors
or blocks of eigenvectors, relying on matrix-vector multiplications rather than
manipulating matrix entries.  
Among these are the Arnoldi iteration and its variants \cite{GoGr06,GoVL96,HWHSG21};
and for symmetric problems, 
the popular Locally Optimal Block Preconditioned Conjugate Gradient (LOBPCG)
\cite{DSYG18,knyazev01}, and the related but more general inverse-free preconditioned 
Krylov subspace methods \cite{GoYe02,QuYe10}.
These methods all use the idea of iteratively projecting the problem onto a Krylov 
subspace of relatively small dimension where dense methods are used to solve a
small eigenvalue problem. 
Additional methods close to this class include the Davidson \cite{davidson75} 
and Jacobi-Davidson \cite{HoNo06} methods from computational chemistry
which use a similar idea, but introduce a preconditioner 
by which the vectors of the
projection subspace are no longer equivalent to a Krylov basis.

An alternate and complementary approach to accelerating eigenvector convergence 
in the power method is based on extrapolation.  The idea is to recombine the latest
update with previous information to form the next iterate in an approximation sequence.  
One of the best known methods in this class is Aitken's acceleration 
\cite[chapter 9]{wilkinson65}, with extensions to vector and $\epsilon$ extrapolation 
methods including \cite{brezinski75,BRZ06,sidi04,sidi08}, to name a few.
Recently, several new methods for accelerating the power method with extrapolation
have been developed, including \cite{BWM20}, in which the power method is recast
as a non-stationary Richardson method; and \cite{NiPo20} which damps the 
largest subdominant eigenmodes to accelerate convergence, and which introduces
the idea of computing a dynamic extrapolation parameter based on a ratio of 
residuals. A similar technique was
used in \cite{PoSc21} to accelerate the Arnoldi iteration.

In  \cite{DSHMRX19}, a power method with an added momentum-type extrapolation 
term was introduced, based on the well known 
heavy ball method of \cite{polyak64}. It was shown 
that this momentum term accelerates the convergence of the power iteration
for positive semidefinite matrices,
and the optimal momentum parameter 
for the acceleration is given by $\beta = \lambda_2^2/4$ where $\lambda_2$ is the 
second largest magnitude eigenvalue of the matrix.  A method to add a beneficial 
momentum term 
without explicit knowledge of $\lambda_2$ was proposed in \cite{DSHMRX19} as the 
Best Heavy Ball method, which relies on multiple matrix-vector multiplications per 
iteration throughout the algorithm. 

To improve upon this method, a delayed momentum power method (DMPower) was proposed in a 
more recent paper \cite{RJRH22}. The method involves a two-phase approach. The first is a 
pre-momentum phase consisting of standard power iterations with inexact deflation, 
at a cost of three matrix-vector multiplies per iteration, to estimate both $\lambda_1$ 
and $\lambda_2$. The second phase runs the method of \cite{DSHMRX19}
with fixed momentum parameter $\beta$ computed with the approximation to $\lambda_2$
from the first phase. 
An analysis is included of how many preliminary iterations are required to obtain
a reliable approximation to $\lambda_2$, based on a priori spectral knowledge.

In this paper, we introduce a dynamic momentum method designed to accelerate the 
power iteration with minimal additional cost per iteration.
In the method proposed herein, the momentum parameter is updated at each iteration 
based on the Rayleigh quotient and  two previous residuals. 
Like the standard power iteration, this method requires only a single matrix-vector 
multiplication per iteration. 
As we will see in section \ref{sec:numerics}, the introduced dynamic method 
outperforms not only the power iteration, but also the static momentum method.
We additionally show in section \ref{sec:inverse} that the method is beneficial
when applied to a shifted inverse iteration. 

We will consider matrix 
$A \in \R^{n\times n}$ with eigenvalues $\lambda_1, \ldots, \lambda_n$ with 
$|\lambda_1| > |\lambda_2| \ge \ldots \ge |\lambda_n|$. The results trivially generalize
to the case where $\lambda_1 = \lambda_2 = \cdots =\lambda_r$ and 
$|\lambda_{r}| > |\lambda_{r+1}| \ge \ldots \ge |\lambda_n|$.
As in \cite{NiPo20}, our proposed method dynamically updates parameters based on the 
detected convergence rate computed by the ratio of 
the last two residuals. 

To fix notation, we can write the power iteration as
\begin{align}\label{eqn:pow}
u_{k+1} = A x_k, \quad x_{k+1} = h_{k+1}^{-1} u_{k+1}, \quad h_{k+1}= \nr{u_{k+1}}.
\end{align}
The momentum method for the power iteration introduced in \cite{DSHMRX19},
takes the form
\begin{align}\label{eqn:mom1}
u_{k+1} = A x_{k} - \beta h_{k}^{-1} x_{k-1},
\quad x_{k+1} = h_{k+1}^{-1} u_{k+1}, \quad h_{k+1}= \nr{u_{k+1}},
\end{align}
where $\beta>0$ is the momentum parameter. As shown in \cite{DSHMRX19} and summarized in section \ref{sec:moment}, an optimal
choice of $\beta$ is $\lambda_2^2 /4$, where it is assumed that 
$|\lambda_2| < |\lambda_1|$. 
Our proposed dynamic method based on iteration \eqref{eqn:mom1} takes the form
\begin{align}\label{eqn:momd}
u_{k+1} = A x_{k} - \beta_k h_{k}^{-1} x_{k-1},
\quad x_{k+1} = h_{k+1}^{-1} u_{k+1}, \quad h_{k+1}= \nr{u_{k+1}}.
\end{align}
This method, described in section \ref{sec:dymo}, 
assigns the parameter $\beta_k$
with minimal additional computation (and no additional matrix-vector multiplies), 
producing a dynamically updated version of 
\eqref{eqn:mom1}.

The remainder of this paper is structured as follows.
Subsections \ref{subsec:prelim}-\ref{subsec:algs} state the basic assumptions
and reference algorithms.
In section \ref{sec:moment} we summarize convergence results for the ``static'' 
momentum method of \cite{DSHMRX19} through the lens of the power iteration 
applied to an augmented matrix.  While this approach was outlined in \cite{DSHMRX19},
our analysis goes a step further, showing that the augmented matrix is defective
under the optimal parameter choice.
In section \ref{sec:dymo} we present the main contributions of this paper: 
our dynamic momentum algorithm 
\ref{alg:dymo}, and an analysis of its convergence and stability.
Numerical results for the method are presented in section \ref{sec:numerics}.
In section \ref{sec:inverse}, we present and discuss 
algorithm \ref{alg:dymoinv} to accelerate the shifted inverse iteration with momentum.

%% -- --------------------------------------------------------------------
\subsection{Preliminaries}\label{subsec:prelim}
%% -- --------------------------------------------------------------------
Our standard assumption throughout the paper is the following.
\begin{assumption}\label{assume:l1diag}
Suppose $A \in \R^{n \times n}$ is diagonalizable and the $n$ eigenvalues of $A$ satisfy
$|\lambda_1| > |\lambda_2| \ge \ldots \ge |\lambda_n|.$
\end{assumption}
Under assumption \ref{assume:l1diag}, let $\{\phi_l\}_{l = 1}^n$ be a set of 
eigenvectors of $A$ so that each $(\lambda_l, \phi_l)$ is an eigenpair of $A$.

In order to analyze the momentum method for $A$, which we will see is equivalent to
a power iteration on an augmented matrix, we will need to make a more general assumption
on the augmented matrix.
\begin{assumption}\label{assume:l1}
Suppose $A \in \R^{n \times n}$ and the $n$ eigenvalues of $A$ satisfy
$|\lambda_1| > |\lambda_2| \ge \ldots \ge |\lambda_n|.$
\end{assumption}
The key difference in assumption \ref{assume:l1} is the matrix is not necessarily
diagonalizable.  In this case we will still refer to the eigenvectors  as 
$\phi_1, \ldots, \phi_n$, but will specify which if any are in fact generalized 
eigenvectors corresponding to a defective eigenspace.

Throughout the paper, $\anorm$ is the Euclidean or $l_2$ norm,
induced by the $l_2$ inner-product denoted by $(\cdot \, , \cdot)$.

% -- --------------------------------------------------------------------
% -- --------------------------------------------------------------------
\subsection{Reference algorithms}\label{subsec:algs}
% -- --------------------------------------------------------------------
Next we state the power iteration \eqref{eqn:pow} and the momentum iteration
\eqref{eqn:mom1} in algorithmic form. The algorithm for the momentum iteration
will require a single preliminary power iteration, and the algorithm for the
dynamic momentum method to be introduced in section \ref{sec:dymo} will require
two preliminary power iterations. 

\begin{algorithm}{Power iteration}
\label{alg:pow}
\begin{algorithmic}
\State{Choose $v_0$, set $h_0 = \nr{v_0}$ and $x_0 = h_0^{-1}v_0$}
\State{Set $v_1 = A v_0$}
\For{$k \ge 0$}
\State{Set $h_{k+1} = \nr{v_{k+1}}$ and  $x_{k+1} = h_{k+1}^{-1} v_{k+1}$}
\State{Set $v_{k+2} = A x_{k+1}$}
\State{Set $\nu_{k+1} = (v_{k+2}, x_{k+1})$ and 
       $d_{k+1} = \|v_{k+2} - \nu_{k+1} x_{k+1}\|$}
\State{{STOP} if $\nr{d_{k+1}} <$ \tt{tol}  }
\EndFor
\end{algorithmic}
\end{algorithm}

The algorithm for the power iteration with momentum assumes knowledge of $\lambda_2$ to assign the parameter $\beta = \lambda_2^2/4$ and implements the iteration \eqref{eqn:mom1}.
\begin{algorithm}{Power iteration with momentum}
\label{alg:hbpow}
\begin{algorithmic}
\State Set $\beta = \lambda_2^2/4$
\State{Do a single iteration of algorithm \ref{alg:pow}} \Comment{$k = 0$}
\For{$k \ge 1$} \Comment{$k \ge 1$}
\State{Set $u_{k+1} = v_{k+1} - (\beta/h_k) x_{k-1}$}
\State{Set $h_{k+1} = \nr{u_{k+1}}$ and $x_{k+1} = h_{k+1}^{-1} u_{k+1}$} 
\State{Set $v_{k+2} = A x_{k+1}$, $\nu_{k+1} = (v_{k+2}, x_{k+1})$ and 
           $d_{k+1} = \nr{v_{k+2} - \nu_{k+1} x_{k+1}}$}
\State{{STOP} if $\nr{d_{k+1}} <$ \tt{tol}  }
\EndFor
\end{algorithmic}
\end{algorithm}

%% -- --------------------------------------------------------------------
\section{Background: the static momentum method}\label{sec:moment}
In this section we will review some results on algorithm \ref{alg:hbpow}, the power iteration 
with momentum. 
To this end, we will also review some standard supporting results on the 
power iteration, algorithm \ref{alg:pow}, in both diagonalizable and defective scenarios. 
These results will be useful to understand each step of the dynamic momentum method.

\subsection{Iteration \eqref{eqn:mom1} as a power iteration with an augmented matrix}\label{subsec:static}
As shown in \cite{DSHMRX19}, the iteration \eqref{eqn:mom1} is equivalent
to the first $n$ rows of the standard power iteration \eqref{eqn:pow} applied to
the augmented matrix
\begin{align}\label{eqn:augsys}
A_\beta = \begin{pmatrix} A & -\beta I \\ I & 0 \end{pmatrix}.
\end{align}

To see this, consider the power iteration on $A_\beta$ starting with $x_0$ in the first 
component (meaning the first $n$ rows) and $y_0$ in the second, then writing
\begin{align}\label{eqn:augiter1}
\begin{pmatrix}u_k\\z_k \end{pmatrix}= 
A_\beta \begin{pmatrix} x_{k-1} \\ y_{k-1} \end{pmatrix} = 
\begin{pmatrix} A x_{k-1} - \beta y_{k-1}\\ x_{k-1} \end{pmatrix}.
\end{align}
Normalizing each component by a scalar $h_k$ (to be discussed below) with $x_k = h_k^{-1} u_k$ and 
$y_k = h_k^{-1} z_k = h_k^{-1} x_{k-1}$ yields the iteration 
\begin{align}\label{eqn:augiter2}
\begin{pmatrix}u_{k+1}\\z_{k+1} \end{pmatrix}= 
A_\beta \begin{pmatrix} x_{k} \\ y_{k} \end{pmatrix} = 
\begin{pmatrix} A x_{k} - \beta y_{k}\\ x_{k} \end{pmatrix}=
\begin{pmatrix} A x_{k} - \beta h_k^{-1}x_{k-1}\\ x_{k} \end{pmatrix}.
\end{align}
The first component in \eqref{eqn:augiter2} agrees with \eqref{eqn:mom1} if we choose $h_k = \norm{u_k}$.
Although this is actually a semi-norm over the tuple $(u_k,z_k)$, it is the most
convenient choice for the sake of computing the Rayleigh quotient corresponding
to the first component at each iteration. 
Hence the equivalence between iteration \eqref{eqn:mom1} and the power iteration 
given by \eqref{eqn:pow} as applied to the augmented matrix \eqref{eqn:augsys}
holds, up to the chosen normalization factor. 
Algorithm \ref{alg:hbpow} explicitly 
performs this iteration starting with $y_0 = 0$ and $\beta = \lambda_2^2/4$, which 
we discuss further below.

The convergence of iteration \ref{eqn:mom1} for general $\beta \in [0, \lambda_1^2/4)$,
$\beta \ne \lambda_i^2/4$, $i = 2, \ldots, n$,
can be quantified in terms of the convergence
of the standard power iteration algorithm \ref{alg:pow}. 
Under assumption \ref{assume:l1diag}, this  
can be summarized as in \cite[Chapter 7]{GoVL96} by
\begin{align}\label{eqn:conv-pow}
\dist\left(\spa \{x_k\},\spa \{\phi_1\} \right) 
= \bigo\left(\left|\f{\lambda_2}{\lambda_1} \right|^k \right), ~\text{ and }~
|\lambda_1 - \nu_k | 
= \bigo\left(\left|\f{\lambda_2}{\lambda_1} \right|^k \right),
\end{align}
which follows by standard arguments from the expansion of initial iterate $u_0$ 
as a linear combination of the $n$ eigenvectors of $A$, namely 
$u_0 = \sum_{l = 1}^n a_l \phi_l, $ by which
\begin{align}\label{eqn:diagpow}
A^k u_0 = a_1 \lambda_1^k\left( \phi_1 + \sum_{l = 2}^n \f{a_l}{a_1} 
\left( \f{\lambda_l}{\lambda_1} \right)^k \phi_l \right).
\end{align}

In the case that assumption \ref{assume:l1} holds and 
$A$ is not diagonalizable, i.e., {\em defective}, 
the power iteration still converges to the dominant eigenpair. 
This is the case for $A_\beta$ when $\beta = \lambda^2/4$ for any subdominant eigenvalue 
$\lambda$ of $A$, as we will see in proposition \ref{prop:spec-aug}. 
For a general defective matrix $A$, if the eigenspace of $\lambda_1$  does not have a full 
set of eigenvectors then the convergence is slow (like $1/k$, 
where $k$ is the iteration count), as shown for instance in \cite[Chapter 9]{wilkinson65}.
If, on the other hand, assumption \ref{assume:l1} holds, $A$ is defective, and
the eigenspace for $\lambda_J$ with $J \ge 2$ lacks a full set of eigenvectors, then
the convergence of algorithm \ref{alg:pow} still agrees with \eqref{eqn:conv-pow}, but
only asymptotically. 
In particular, from \cite[Chapter 9]{wilkinson65}, 
if for $|\lambda_J| < |\lambda_1|$ we have
$\lambda_J = \lambda_{J+1}$ and the corresponding
eigenspace has geometric multiplicity 1, then letting $\phi_{J+1}$ be a generalized 
eigenvector with $A \phi_{J+1} = \lambda_J \phi_{J+1} + \phi_J$,
in place of  \eqref{eqn:diagpow} we have
\begin{align}\label{eqn:defecpow}
A^k u_0 = a_1 \lambda_1^k\left( \phi_1 +  
\f{a_{J+1}}{a_{1}}\left( \f{k \lambda_J^{k-1}}{\lambda_1^k} \right) \phi_J +  
\sum_{l = 2}^n \f{a_l}{a_1} 
\left( \f{\lambda_l}{\lambda_1} \right)^k \phi_l \right).
\end{align}
Noting that $(k \lambda_J^{k-1}/\lambda_1^k) / 
((k-1) \lambda_J^{k-2}/\lambda_1^{k-1})  \rightarrow \lambda_J/\lambda_1 $
as $k \goto \infty$, we have the same asymptotic convergence rate as in the 
non-defective case.
This is important for the analysis of algorithm \ref{alg:hbpow}
since as shown in the next proposition, the augmented matrix $A_\beta$ is defective
whenever $\beta = \lambda^2/4$ for any eigenvalue $\lambda$ of $A$. 

\subsection{Spectrum of the augmented matrix}\label{subsec:augspec}
By the equivalence between the first component
of the power iteration on $A_\beta$ and algorithm \ref{alg:hbpow} as shown in
\eqref{eqn:augiter1}-\eqref{eqn:augiter2}, 
the convergence rate of the momentum accelerated method of iteration \eqref{eqn:mom1} 
depends on ratio of the two largest magnitude eigenvalues of $A_\beta$.
In order to understand the convergence properties of algorithm \ref{alg:hbpow} 
and later our dynamic version of this method, 
the following proposition describes the spectral decomposition  of 
$A_\beta$ in terms of the eigenvalues and eigenvectors of $A$.

\begin{prop}\label{prop:spec-aug}
Suppose $A$ satisfies assumption \ref{assume:l1diag}.
Then the $2n$ (counting multiplicity) eigenvalues of $A_\beta$ are given by
\begin{align}\label{eqn:mulambda}
\mu_{\lambda_\pm} = \f 1 2 \left( \lambda \pm \sqrt{\lambda^2 - 4\beta} \right),
\quad \lambda \in \{\lambda_1, \ldots, \lambda_n\}.
\end{align}

In the case that $\lambda^2 - 4 \beta \ne 0$, the eigenvectors of $A_\beta$ 
corresponding to each eigenvalue $\mu = \mu_{\lambda_{\pm}}$ are given by
\begin{align}\label{eqn:aug-vect}
\psi_{\lambda_\pm} = \begin{pmatrix} (\mu_{\lambda_{\pm}}) \phi \\ \phi\end{pmatrix},
\end{align}
where $\phi$ is the eigenvector of $A$ corresponding to eigenvalue $\lambda$.

In the case that $\beta =\lambda^2/4 > 0$, the matrix $A_\beta$ is not diagonalizable.
Moreover, if $\lambda $ is an eigenvalue of multiplicity $m$ of $A$, then the eigenvalue 
$\mu_\lambda = \lambda/2$ of $A_\beta$ has algebraic multiplicity $2m$ and geometric 
multiplicity $m$.
\end{prop}

We restrict our attention to $\beta > 0$ as iteration \eqref{eqn:mom1} reduces to
\eqref{eqn:pow} if $\beta = 0$.
Before the proof of proposition~\ref{prop:spec-aug}, we include a corollary that follows immediately from its conclusions.
\begin{cor}\label{cor:abeta}
If $A$ satisfies assumption \ref{assume:l1diag} and $\beta \in (0, \lambda_1^2/4)$, 
then $A_\beta$ as given by \eqref{eqn:augsys} satisfies assumption  \ref{assume:l1}.
\end{cor}

Together, for symmetric matrices,
proposition \ref{prop:spec-aug} and corollary \ref{cor:abeta}
show that as the power iteration applied to the 
augmented matrix $A_\beta$ converges, the first component of the eigenvector converges 
to the dominant eigenvector of 
$A$ for any $\beta \in (0,\lambda_1^2/4)$. 
If $\beta = \lambda^2/4$ for any nonzero $\lambda = \lambda_2, \ldots, \lambda_n$, 
then the matrix $A_\beta$ is defective, but courtesy of \eqref{eqn:defecpow}, 
the power iteration will converge asymptotically at the same rate as in the diagonalizable 
case as given by \eqref{eqn:diagpow}, applied to the eigenvalues of $A_\beta$.

\begin{proof}
The eigenvectors of $A_\beta$ are related to the eigenvectors of $A$ 
by noting that if $\phi$ is an eigenvector of $A$ with eigenvalue $\lambda$ then solving 
\[
A_\beta \begin{pmatrix} \mu \phi \\ \phi \end{pmatrix} 
= \mu\begin{pmatrix} \mu \phi \\ \phi\end{pmatrix}, \quad ~\text{ which reduces to }~
\begin{pmatrix}(\mu\lambda - \beta)\phi \\ \mu \phi \end{pmatrix}
= \mu\begin{pmatrix} \mu \phi \\ \phi\end{pmatrix},
\]
for $\mu \in \C$, yields the quadratic equation $\mu^2 - \lambda \mu + \beta = 0$. If $\beta \neq \lambda^2/4$,  the $2n$ eigenvalues of $A_\beta$ are given by \eqref{eqn:mulambda}, and the corresponding eigenvectors are given by \eqref{eqn:aug-vect}. 

On the other hand, if $\beta = \lambda^2/4$ where $\lambda$ is an eigenvalue of $A$ with algebraic multiplicity 1, then the quadratic equation
$\mu^2 - \lambda \mu + \beta = 0$ has a repeated root $\mu = \lambda/2$. To find the eigenvector(s) associated with $\mu$, 
we can express the equation for null-vectors of $A_\beta - \mu I$ as 
\[
\begin{pmatrix}A - \f{\lambda}{2} I & -\f{\lambda^2}{4} I 
            \\ I                 & - \f{\lambda}{2} I \end{pmatrix}
\begin{pmatrix} u \\ v\end{pmatrix} = 
\begin{pmatrix} 0 \\ 0\end{pmatrix}.
\]
From the second component of the equation, $u = \f{\lambda}{2} v$. Applying this to the
first component yields $(A - \f{\lambda}{2}I) \f{\lambda}{2}v - \f{\lambda^2}{4}v = 0$,
or $Av = \lambda v$. This implies that $v$ must be an eigenvector of $A$ corresponding to 
eigenvalue $\lambda$. Therefore, the eigenspace for $A_\beta$ corresponding
to the repeated eigenvalue $\mu = \lambda/2$ has dimension $1$.

More generally, if $\lambda$ is an eigenvalue of algebraic and geometric multiplicity $m>1$,
then the argument above can be applied to each eigenpair $(\lambda,\hat \phi_i)$, 
$i = 1, \ldots, m$, where $\{\hat \phi\}_{i = 1}^m$ is some basis for the eigenspace
corresponding to $\lambda$.
Then for $\beta = \lambda^2/4$, $A_\beta$ has an eigenvalue 
$\mu = \lambda/2$ with algebraic multiplicity $2m$ but with geometric multiplicity $m$. 
\end{proof}

From \eqref{eqn:mulambda} of proposition \ref{prop:spec-aug} 
we have three cases for each pair of %the $2n$ 
eigenvalues of $A_\beta$ corresponding 
to a real eigenvalue of $A$, determined by 
the sign of the discriminant in
\eqref{eqn:mulambda}.
Define $\mu_\lambda$ as the larger magnitude eigenvalue
of $A_\beta$ corresponding to eigenvalue $\lambda$ of $A$, and $\widehat \mu_\lambda$
as the smaller magnitude corresponding eigenvalue, in the case
that $\mu_{\lambda_\pm}$ are real. 
If $\mu_{\lambda_{\pm}}$ are complex, define $\mu_\lambda$ as having the positive 
imaginary component.
Then
\begin{align}\label{eqn:mucase1}
(\lambda/2)^2 &\ge \beta: &&  \mu_\lambda = 
\f 1 2 \left( \lambda +\sign(\lambda)  \sqrt{\lambda^2 -4 \beta} \right)
\\ \label{eqn:mucase2}
(\lambda/2)^2 &= \beta: &&  \mu_\lambda = \f 1 2 \lambda
\\ \label{eqn:mucase3}
(\lambda/2)^2 &\le \beta: &&  \mu_\lambda = {\sqrt{\beta}}e^{\imath \theta}, ~\text{ with }~
\theta = \arctan \left(\sqrt{\f{4\beta}{\lambda^2} -1} \right),
\end{align}
where \eqref{eqn:mucase2} agrees with both \eqref{eqn:mucase1} and \eqref{eqn:mucase3}
at $\beta = (\lambda/2)^2$, and is separately enumerated only for emphasis. 
In \eqref{eqn:mucase3}, it is understood that $\theta = \pi/2$ when $\lambda =0$.  Based on \eqref{eqn:mucase3}, 
we see $\beta \ge \lambda_1^2/4$ causes all real
eigenvalues of $A_\beta$ to have
equal magnitude $\sqrt{\beta}$.
If $A$ has complex eigenvalues, the complete set of eigenvalues can still be 
given by 
$\f 1 2 \left( \lambda \pm  \sqrt{\lambda^2 -4 \beta} \right)$,
applied to each eigenvalue $\lambda$ of $A$, however the quantity in the square
root may be complex.

We can now summarize the convergence properties of the standard power iteration 
\eqref{eqn:pow} applied to the augmented matrix $A_\beta$ given by
\eqref{eqn:augsys} hence iteration \eqref{eqn:mom1} 
for symmetric matrices $A$ as follows.
An alternate approach based on Chebyshev polynomials 
shown for positive semidefinite matrices 
can be found in \cite{DSHMRX19}.

\begin{cor}\label{cor:conv-aug}
For $0 < \beta < \lambda_1^2/4$,
the power iteration \eqref{eqn:pow} implemented in algorithm \ref{alg:pow} 
applied to the augmented matrix $A_\beta$
of \eqref{eqn:augsys} for symmetric matrix $A$ converges at the rate 
\begin{align}\label{eqn:cor-augrate}
\f{|\mu_{\lambda_2}|}{|\mu_{\lambda_1}|} =
\left\{ \begin{array}{ll}
 \f{2\sqrt{\beta}}
   {|\lambda_1| +  \sqrt{\lambda_1^2 - 4\beta}}, 
& \lambda_2^2/4 < \beta < \lambda_1^2/4 \\ \\
 \f{ |\lambda_2| +  \sqrt{\lambda_2^2 - 4\beta} }
    {|\lambda_1| +  \sqrt{\lambda_1^2 - 4\beta} }, 
& 0 \le \beta < \lambda_2^2/4,
\end{array}\right. 
\end{align}
and asymptotically at the rate
\begin{align}\label{eqn:cor-momrate}
\f{|\mu_{\lambda_2}|}{|\mu_{\lambda_1}|} \goto 
 \f{ |\lambda_2| }{|\lambda_1| + \sqrt{\lambda_1^2 - 4\beta} }
   = \f {r}{1 + \sqrt{1 - r^2}}, ~\text{ with }~ 
r = |\lambda_2 /\lambda_1| , ~\text{ for }~ \beta = \lambda_2^2/4.
\end{align}
The choice of $\beta$ that optimizes the asymptotic convergence rate is 
$\beta = \lambda_2^2/4$,
for which the power iteration applied to $A_\beta$ and the power iteration with 
momentum algorithm \ref{alg:hbpow} applied to $A$ 
converge asymptotically at the rate given 
by \eqref{eqn:cor-momrate}.
\end{cor}
\begin{figure}
\includegraphics[trim = 0pt 0pt 0pt 0pt,clip = true, width = 0.4\textwidth]
{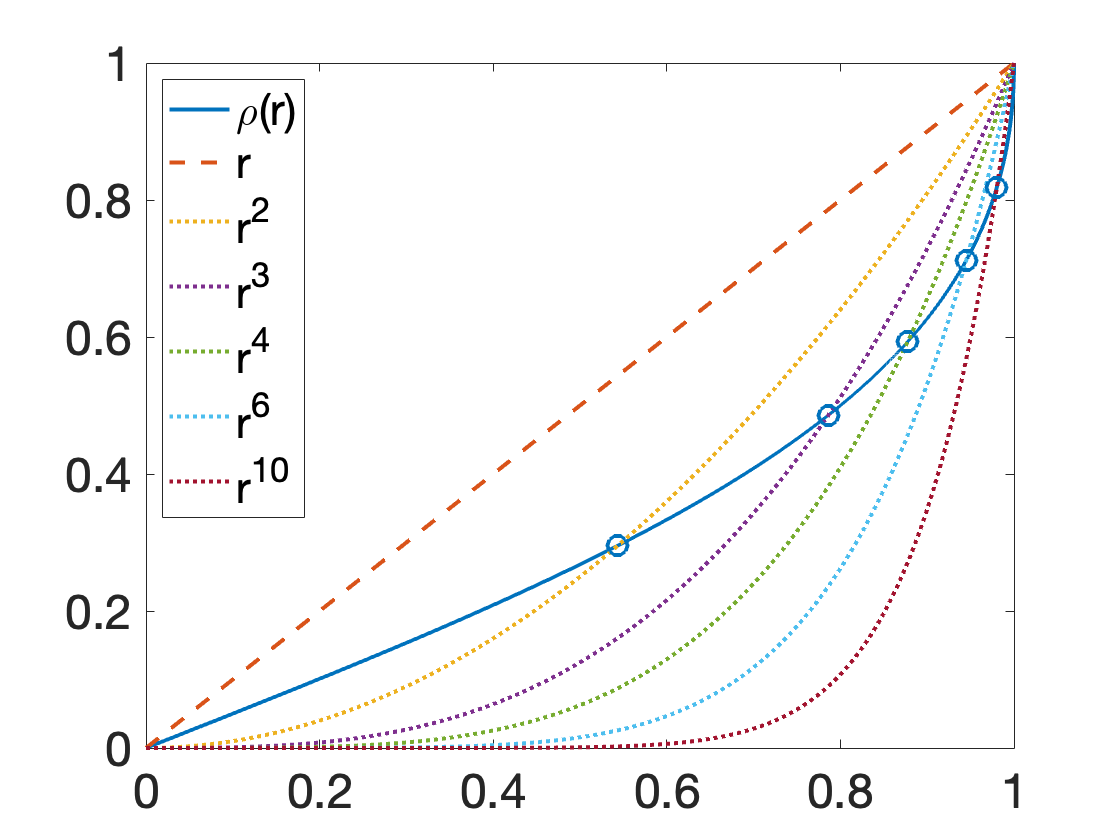}~\hfil~
\includegraphics[trim = 0pt 0pt 0pt 0pt,clip = true, width = 0.4\textwidth]
{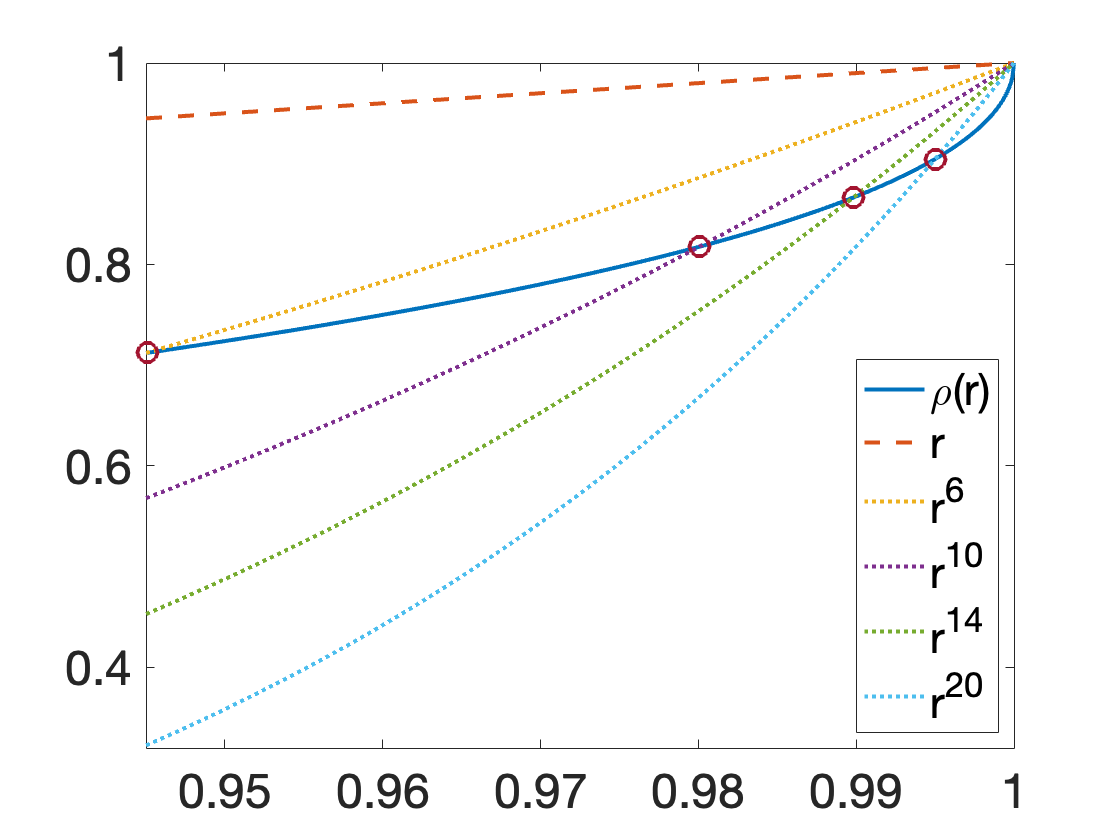}
\caption{A comparison of $\rho(r)$ vs. $r^p$ for $\rho(r) = r/(1+\sqrt{1-r^2})$, 
the rate given in \eqref{eqn:cor-momrate}.
Left: $\rho(r)$ compared with $r^p$, for $p = 1,2,3,4,6,10$. 
Right: a detail plot of $\rho(r)$ compared with $r^p$, for $p = 6,10,14,20$.
The crossings between $\rho(r)$ and $r^p$ are marked in each plot.
\label{fig:rhoplot}}
\end{figure}
As visualized in figure \ref{fig:rhoplot},
the rate given by \eqref{eqn:cor-momrate} is less than $r$ for $r \in (0,1)$, 
less than $r^3$ for $r \in (0.786,1)$,
less than $r^4$ for $r \in (0.878,1)$, 
and less than $r^6$ for $r \in (0.945,1)$, etc.
Hence the smaller the spectral gap in $A$, namely the closer $r = |\lambda_2/\lambda_1|$ is 
to 1,  the more beneficial it is to apply the acceleration.

\begin{remark}
Corollary \ref{cor:conv-aug} shows that the the power iteration applied to the
augmented matrix $A_\beta$ of \eqref{eqn:augsys} converges to the dominant eigenpair
$(\mu_\lambda,\psi_\lambda)$ of $A_\beta$ at a faster rate then the power iteration 
applied to matrix $A$ converges to its dominant eigenpair $(\lambda, \phi)$. 
Proposition \eqref{prop:spec-aug} shows that the dominant
eigenvector of $A$ is the first component (the first $n$ entries) of the dominant
eigenvector of $A_\beta$ for symmetric $A$. As the momentum method \eqref{eqn:mom1} 
generates the first component of the power iteration for $A_\beta$ (using a different
normalization factor), this method approximates the dominant eigenvector of $A$, and 
converges at the rate described in corollary \ref{cor:conv-aug}. The dominant 
eigenvalue $\lambda$ of $A$ can then be recovered by taking a Rayleigh quotient with 
the approximate eigenvector. In practice the augmented matrix $A_\beta$ is never formed;
it is used here as a tool in the analysis of iteration \eqref{eqn:mom1}.
\end{remark}

\begin{proof}
The main technicality in the proof of \eqref{eqn:cor-augrate} is verifying that
$|\widehat \mu_{\lambda_1}| < |\mu_{\lambda_2}|$. Then 
from standard theory, e.g., \cite{GoVL96}, the (asymptotic) rate of convergence
to the eigenvector $\psi_1$ corresponding to $\mu_{\lambda_1}$ is given by
$|\mu_{\lambda_2}/\mu_{\lambda_1}|$.

Without loss of generality, 
suppose $\lambda_1 > 0$. Then for any $\beta \in (0,\lambda_1^2/4)$, we have 
\[
\widehat{\mu}_{\lambda_1} = 
\frac{1}{2}\left(\lambda_1 - \sqrt{\lambda_1^2 - 4\beta}\right).
\] 
By \eqref{eqn:mucase1}-\eqref{eqn:mucase3}, we have $|\mu_{\lambda_2}| \ge \sqrt{\beta}$.
Hence to see that $|\widehat \mu_{\lambda_1}| < |\mu_{\lambda_2}|$, it suffices to show
that $|\widehat \mu_{\lambda_1}| \le \sqrt{\beta}$. This is true since
\[
 \sqrt{\beta} - \widehat{\mu}_{\lambda_1} 
         = \sqrt{\beta} - \f{\lambda_1}{2} + \sqrt{\f{\lambda_1^2}{4} - \beta}
         = \sqrt{\frac{\lambda_1}{2} - \sqrt{\beta}}\left(\sqrt{\f{\lambda_1} {2} +\sqrt{\beta}} - \sqrt{\frac{\lambda_1}{2} -\sqrt{\beta}}\right) \ge 0.
\]
The result \eqref{eqn:cor-augrate} then follows directly from 
\eqref{eqn:mucase1}-\eqref{eqn:mucase3}. 

Next we show the asymptotic optimality of $\beta = \lambda_2^2/4$.  
For this purpose, we consider the convergence rate ~\eqref{eqn:cor-augrate} 
as a function of $\beta$ (for $\beta \ne \lambda_2^2/4$) defined as: 
\[
h(\beta) =
\left\{ \begin{array}{ll}
	\f{2\sqrt{\beta}}
	{|\lambda_1| + \sqrt{\lambda_1^2 - 4\beta}}, 
	& \lambda_2^2/4 < \beta < \lambda_1^2/4 \\ \\
	\f{ |\lambda_2| +  \sqrt{\lambda_2^2 - 4\beta} }
	{|\lambda_1| + \sqrt{\lambda_1^2 - 4\beta} }, 
	& 0 < \beta \le \lambda_2^2/4.
\end{array}\right. 
\]
By direct calculation, we get $h'(\beta) >0$ for $\beta \in (\lambda_2^2/4, \lambda_1^2/4)$, 
that is, the convergence rate is increasing with respect to $\beta$.  
For $\beta \in (0,\lambda_2^2/4)$, we have $h'(\beta) < 0$, so the convergence rate is 
decreasing on $\beta$. Hence by continuity, 
$h(\beta)$ achieves a minimum at  $\beta = \lambda_2^2/4$. 
We note that $h(\beta)$ agrees with the convergence rate $|\mu_{\lambda_2}/\mu_{\lambda_1}|$
except when $\beta = \lambda_2^2/4$. When $\beta = \lambda_2^2/4$, the agreement is
only asymptotic, that is  $|\mu_{\lambda_2}/\mu_{\lambda_1}| \goto h(\beta)$.
\end{proof}

\begin{figure}
\includegraphics[trim = 0pt 0pt 0pt 0pt,clip = true, width = 0.48\textwidth]{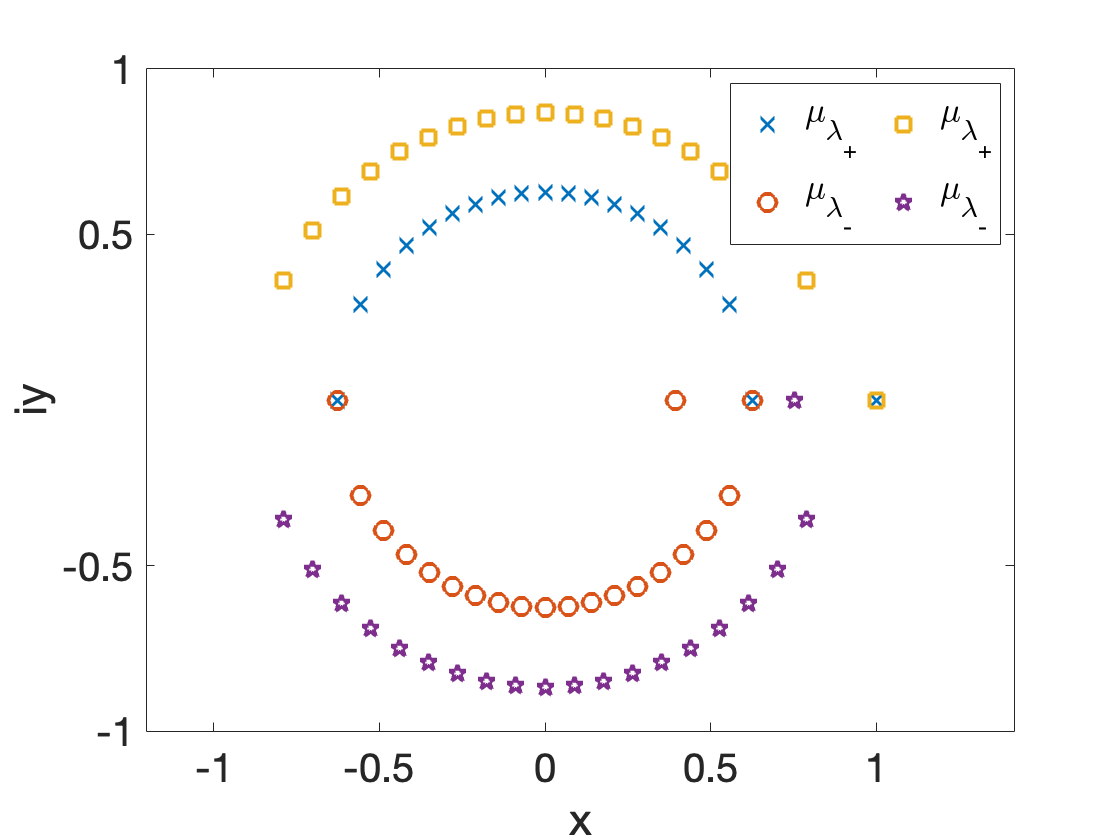}~
\includegraphics[trim = 0pt 0pt 0pt 0pt,clip = true, width = 0.48\textwidth]{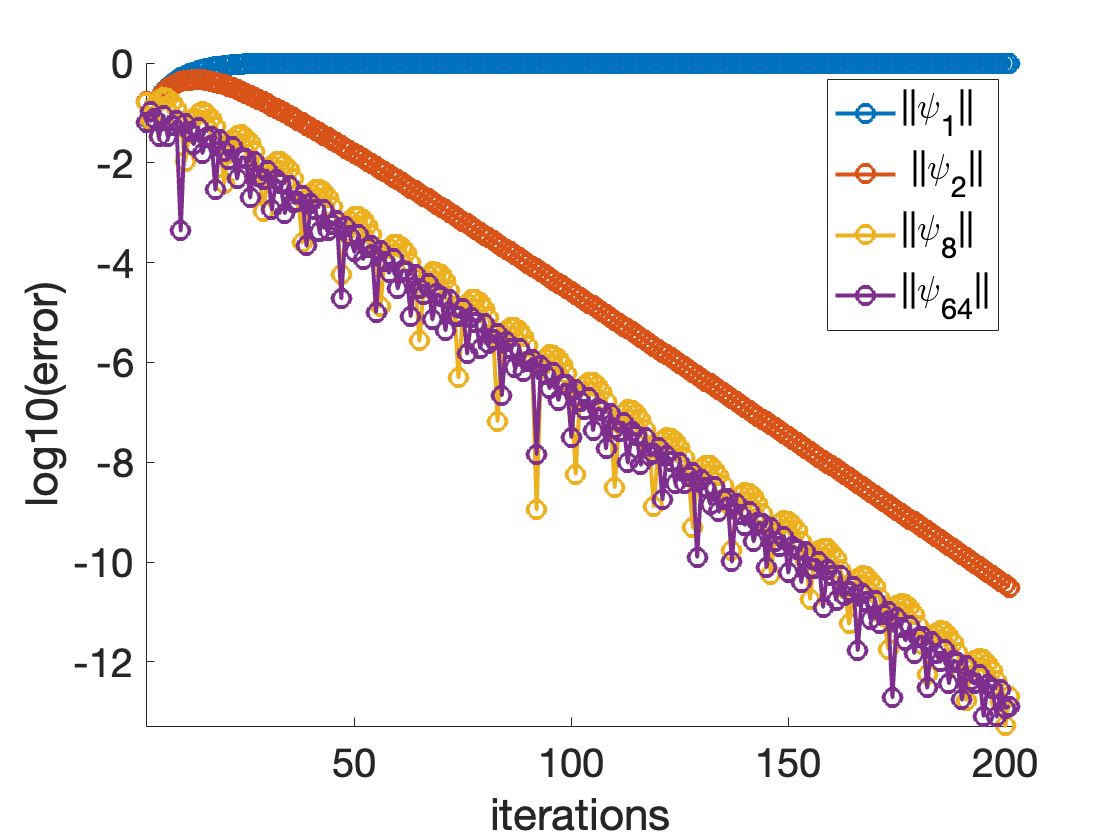}
\caption{Left: The ratio of eigenvalues $\mu_{\lambda_+}/\mu_{\lambda_1}$ and 
 $\mu_{\lambda_-}/\mu_{\lambda_1}$ of the augmented matrix
$A_\beta$ for 
$A = \diag(10:-1:-9)$ with $\beta = 9^2/4$ (inner circle) and 
$\beta = 9.9^2/4$ (outer circle).
Right: convergence of the eigenmodes $\psi_1, \psi_2, \psi_8$ and $\psi_{64}$ of the
augmented matrix $A_\beta$ for $A = \diag(100:-1:1)$ and $\beta = 99^2/4$.  All the 
subdominant modes converge at the same rate, but with increasing oscillation. 
\label{fig:mupics}
}
\end{figure}

Two interesting observations follow from this analysis.
First, as shown in section \ref{sec:numerics}, as well as in the numerical results of 
\cite{DSHMRX19}, iteration \eqref{eqn:momd} with a well chosen dynamically assigned 
sequence of parameters $\beta_k$, for which in general 
$\beta_k \ne \lambda_2^2/4$, can converge faster than the iteration \eqref{eqn:mom1} with the
optimal parameter $\beta = \lambda_2^2/4$.  
This can be explained by the above analysis which shows the 
optimal parameter is only asymptotically optimal. Our results of sections
\ref{sec:numerics} and \ref{sec:inverse} show that a close but inexact 
approximation to this parameter can give a better rate of convergence, at least in the
preasymptotic regime.

Second, for $\beta \in [\lambda_2^2/4, \lambda_1^2/4)$, except for $\mu_{\lambda_1}$ and $\widehat{\mu}_{\lambda_1}$ all the remaining $2n-2$ (complex) eigenvalues of $A_\beta$ (corresponding to the eigenvalues 
$\lambda_2, \ldots, \lambda_n$ of $A$) have the same magnitude $\sqrt{\beta}$ according to \eqref{eqn:mucase3}. However, as the corresponding
eigenvalues $\lambda_j$ of $A$ with $|\lambda_j| < |\lambda_2|$ 
decrease in magnitude, the argument $\theta$ in \eqref{eqn:mucase3} increases.
This causes oscillatory convergence at an increasing rate of oscillation for the
subdominant modes.
This is illustrated in figure \ref{fig:mupics}: the left plot shows the ratio of 
eigenvalues $\mu_{\lambda_{\pm}}/\mu_{\lambda_1}$ 
of $A_\beta$ plotted on the complex plane for $\beta = 9^2/4$ (inner circle) and
$\beta = 9.9^2/4$ (outer circle), where $A = \diag(10:-1:-9)$.  The right plot
shows the magnitude of the first, second, eighth and 64th eigenmodes of the power 
iteration algorithm \ref{alg:pow} applied to the augmented matrix $A_\beta$ for
$A = \diag(100:-1:1)$ with $\beta = 99^2/4$. The plots agree with the above analysis: the
modes all decay at the same rate, but the modes of $A_\beta$ corresponding to the 
eigenmodes of $A$ with smaller magnitude eigenvalues have larger imaginary parts, and 
their convergence is more oscillatory. 
The above analysis also shows that if $\beta\ge \lambda_1^2/4$, then all eigenvalues 
of $A_\beta$ have the same magnitude. 
Therefore, if $\beta\ge \lambda_1^2/4$, the augmented matrix $A_\beta$ does not satisfy 
assumption~\ref{assume:l1}, and neither the power iteration applied to $A_\beta$, nor
iteration \ref{eqn:mom1} applied to $A$, will converge. 

%% -- --------------------------------------------------------------------
\section{Dynamic momentum method}\label{sec:dymo}
We would like to use the acceleration of the momentum algorithm \ref{alg:hbpow}, but
without the a priori knowledge of $\lambda_2$.  A method for determining an effective
sequence of momentum parameters is presented in 
\cite[algorithm 3]{DSHMRX19}, called the Best Heavy Ball method. This method however 
requires five matrix-vector multiplications per iteration, as compared to the single 
matrix-vector multiplication per iteration required by the standard power iteration 
algorithm \ref{alg:pow} or the momentum accelerated power iteration \cite{DSHMRX19}
presented here as algorithm \ref{alg:hbpow}. 
This is improved upon in the DMPower algorithm of \cite{RJRH22} which uses inexact
deflation \cite[chapter 4]{saad2011} in a preliminary iteration to approximate $\lambda_2$. However,
the method is sensitive to the approximation of $\lambda_2$, and ensuring the 
approximation is good enough again requires a priori knowledge of the spectrum.
Additionally, the preliminary iteration is more computationally expensive, 
requiring 3 matrix-vector multiplications per iteration.

Our approach for approximating the momentum 
parameter $\beta = \lambda_2^2/4$ 
does not require any additional matrix-vector multiplication per iteration.
We obtain an expression for $r_{k+1}$, as an approximation of $r = |\lambda_2/\lambda_1|$ from the detected residual convergence rate $\rho_k = d_{k+1}/d_k$
by inverting the optimal convergence rate \eqref{eqn:cor-momrate} for 
$r$ in terms of $\rho$.
This is justified in lemma \ref{lem:rstab}.
We then approximate $\lambda_2$ by $r_{k+1}$ multiplied by 
the (computed) Rayleigh quotient approximation to $\lambda_1$, which yields the approximated momentum parameter $\beta_k$.  The resulting dynamic momentum algorithm is presented below.

%% -- --------------------------------------------------------------------
%% -- --------------------------------------------------------------------
\begin{algorithm}{Dynamic momentum}
\label{alg:dymo}
\begin{algorithmic} 
\State {Do two iterations of algorithm \ref{alg:pow}} \Comment{$k = 0,1$}
\State{Set $r_{2} = \min\{d_2/d_{1},1\}$}
\For{$k \ge 2$} \Comment{$k \ge 2$}
\State{Set $\beta_k = \nu_{k}^2 r_{k}^2/4$}
\State{Set $u_{k+1} = v_{k+1} - (\beta_k/h_k) x_{k-1}$}
\State{Set $h_{k+1} = \nr{u_{k+1}}$ and $x_{k+1} = h_{k+1}^{-1} u_{k+1}$}
\State{Set $v_{k+2} = A x_{k+1}$, $\nu_{k+1} = (v_{k+2}, x_{k+1})$ and
           $d_{k+1} = \nr{v_{k+2} - \nu_{k+1} x_{k+1}}$}
\State{Update $\rho_k = \min\{d_{k+1}/d_k,1\}$  and
              $r_{k+1} = 2\rho_k/(1 + \rho_k^2)$}
\State{{STOP} if $\nr{d_{k+1}} <$ \tt{tol}  }
\EndFor
\end{algorithmic}
\end{algorithm}
%% -- --------------------------------------------------------------------
%% -- --------------------------------------------------------------------

%% -- --------------------------------------------------------------------
%% -- --------------------------------------------------------------------

Lemma \ref{lem:rstab} and remark \ref{rem:betastab} 
in the next section show that assigning $r_{k+1}$ by 
$r_{k+1} = {2\rho_k}/(1+\rho_k^2)$, obtained by inverting the 
asymptotic convergence rate \eqref{eqn:cor-momrate} of the optimal parameter $\beta$, 
gives a stable approximation to $r$ and hence to $\beta$.
In fact, the approximation becomes increasingly stable as $r$ gets closer to unity.

The next remark describes the role of the subdominant eigenmodes in the residual.
%% -- --------------------------------------------------------------------

\begin{remark}\label{rem:resi}
The residual $d_k$ as given in algorithms \ref{alg:pow}, \ref{alg:hbpow} and 
\ref{alg:dymo} is given by
$d_k = \nr{A x_{k} - \nu_k x_k}$ where the Rayleigh quotient $\nu_k$ is given by
$(Ax_k,x_k)$. Let $x_k = \sum_{l = 1}^n \alpha_l^{(k)} \phi_l$ where $\{\phi_l\}_{l=1}^n$ is the eigenbasis of $A$. Then 
\begin{align}\label{eqn:resi} 
d_k = \nr{\sum_{l = 1}^n (\lambda_l \alpha_l^{(k)} \phi_l)
         -\sum_{l = 1}^n (\nu_k \alpha_l^{(k)} \phi_l)  } 
    = \nr{\sum_{l = 1}^n (\lambda_l - \nu_k) \alpha_l^{(k)} \phi_l}.
\end{align}
The detected convergence rate $\rho_k$ is given by
\begin{align}\label{eqn:rhok-pa}
\rho_k = \f{d_{k+1}}{d_k} = 
\f{\nr{\sum_{l = 1}^n (\lambda_l - \nu_{k+1}) \alpha_l^{(k+1)} \phi_l}}
{\nr{\sum_{l = 1}^n (\lambda_l - \nu_k) \alpha_l^{(k)} \phi_l}}.
\end{align}
We will consider the preasymptotic regime to be that in which $\lambda_1 - \nu_k$ 
is not negligible in comparison to the coefficients $\alpha_l^{(k)}, ~l > 1$, 
which will be seen to decay.
In the asymptotic regime, we have $\nu_k \approx \lambda_1$ hence 
\eqref{eqn:rhok-pa} reduces for practical purposes to
\begin{align}\label{eqn:rhoka}
\rho_k \approx 
\f{\nr{\sum_{l = 2}^n (\lambda_l - \nu_{k+1}) \alpha_l^{(k+1)} \phi_l}}
{\nr{\sum_{l = 2}^n (\lambda_l - \nu_k) \alpha_l^{(k)} \phi_l}}.
\end{align}

In the usual analysis of the power iteration, coefficients $\alpha_l^{(k+1)}$ 
decay like $\lambda_l/\lambda_1$ at each iteration as in \eqref{eqn:diagpow}, 
hence eventually \eqref{eqn:rhoka} is dominated by the maximal such ratio $\lambda_2/\lambda_1$.  
In contrast, in the case of the augmented matrix $A_{\beta_k}$, 
for each of the eigenmodes with $\lambda^2/4 < \beta_k$, each of the corresponding
eigenvalues has the same magnitude; and, as shown in \eqref{eqn:mucase3}
increasing imaginary parts as corresponding eigenvalues of $A$ decrease.   Hence it
is not necessarily the case that the second eigenmode will dominate \eqref{eqn:rhoka}
through most of the iteration.
The oscillation of the subdominant modes
is the main reason we will see the sequence of convergence rates $\rho_k$
fluctuate in the dynamic algorithm.

However, the stability 
of $r_k$ with respect to $\rho_k$ shown in lemma \ref{lem:rstab} controls the 
oscillations in $r_k$ with respect to $\rho_k$, and substantially damps them in the 
case that $r = |\lambda_2/\lambda_1|$ is close to unity. 
In this case there is a more substantial relative gap between the convergence rate for 
the second eigenmode and the higher frequency modes, so long as some of the $\beta_k$
satisfy $\beta_k < \lambda_2^2/4$, which is generally the case.
Then the second eigenmode does (eventually) tend to dominate the residual.
A further discussion of the coefficients $\alpha_l^{(k)}$ will be given in remark \ref{rem:resicoef}, 
where it will be shown that $\alpha_l^{(k)}$ is controlled by the product of 
eigenvalues of the sequence of augmented matrices $A_{\beta_k}$ corresponding 
to $\lambda_l$ of $A$.
The differences in convergence behavior between
smaller and larger values of $r$ are highlighted in section \ref{sec:inverse}.
\end{remark}

The following subsection takes into account the nontrivial detail that the dynamic algorithm
\ref{alg:dymo} differs from a standard power method in that a different augmented matrix
$A_{\beta_k}$ is applied at each iteration.

%% -- --------------------------------------------------------------------
\subsection{Convergence theory}\label{subsec:momconv}
%% -- --------------------------------------------------------------------
In Section~\ref{subsec:static}, we interpret the convergence of the momentum method with 
constant $\beta$ as a power method applied to the augmented matrix $A_\beta$. However, this perspective no longer precisely holds for algorithm \ref{alg:dymo} as the 
parameter $\beta_k$ is subject to change at each step. Consequently, the corresponding 
augmented matrix $A_{\beta_k}$ changes at each step as well. 
This presents a significant challenge in the analysis of the dynamic momentum algorithm. 

For ease of presentation, we next define some notation to be used throughout the remainder
of this section. Let 
\[
A^{(0)} = \begin{pmatrix} A & 0 \\ I & 0 \end{pmatrix}, \quad \text{ and }\quad 
A^{(j)} = A_{\beta_j} = \begin{pmatrix} A & -\beta_j I \\ I & 0 \end{pmatrix}, ~ j \ge 1,
\]
where $A^{(0)}$ is the augmented matrix with $\beta = 0$.
As in subsection \ref{subsec:prelim}, let $\{\phi_l\}_{l=1}^n$ be an eigenbasis of $A$, 
with corresponding eigenvalues $\{ \lambda_l\}_{l=1}^n$. 
For each eigenpair $(\lambda_l, \phi_l)$ of $A$, denote 
$(\mu^{(j)}_l, \psi^{(j)}_l)$ the corresponding eigenpair of $A^{(j)}$ where 
$\mu_l^{(j)}$ is the eigenvalue with larger magnitude defined in 
\eqref{eqn:mucase1}-\eqref{eqn:mucase3}. 
Then by \eqref{eqn:aug-vect} 
\[
\psi^{(j)}_l = \begin{pmatrix} \mu^{(j)} \phi_l \\ \phi_l \end{pmatrix}, 
\] 
for $j\ge 0$ with $\mu_l^{(0)} = \lambda_l$.

In the first technical lemma of this section we show the effect of applying a sequence of
augmented matrices with changing parameter $\beta_j$ to each eigenmode of $A$.
\begin{lemma}\label{lem:dm-by-mode}
Let $A$ satisfy assumption \ref{assume:l1diag}, let $(\lambda, \phi)$ be an eigenpair
of $A$,  
and let $\mu^{(j)}$ be the corresponding eigenvalue of $A^{(j)}$, 
as in proposition \ref{prop:spec-aug}.
Let $\delta_{ik} = \mu^{(i)} - \mu^{(k)}$. 
Define $\cP^{i}(\mu)$ to be a product of $i$ terms $\mu^{(k)}$, where $1 \le k \le j$, 
and $\cP^i(\delta)$ to be a product of $i$ terms $\delta_{kp}$, where $0 \le k,p \le j$.
Then
\begin{align}\label{eqn:lem-by-mode}
A^{(j)}\cdots A^{(0)} \begin{pmatrix} \phi \\ 0 \end{pmatrix}
= \left( \prod_{i = 1}^j \mu^{(i)}  +
   \sum_{k=1}^{j-1} \delta_{k-1,k}\prod_{i = 1, i \ne k}^j \mu^{(i)} 
   + \sum_{i = 2}^{j-1} \cP^i(\delta)\cP^{j-i}(\mu)\right)
\begin{pmatrix} \mu^{(j)} \phi \\ \phi \end{pmatrix}
\nonumber \\
+ \left( \delta_{j-1,j} \prod_{i=1}^{j-1}\mu^{(i)} + \sum_{i = 2}^j 
\cP^i(\delta) \cP^{j-i}(\mu) \right)
\begin{pmatrix} \mu^{(0)} \phi \\ \phi \end{pmatrix}.
\end{align}
\end{lemma}
This lemma shows that applying the sequence of augmented matrices $A^{(j)}\cdots A^{(0)}$
to each eigenmode of $A$ yields 
a perturbation to multiplying the 
eigenmode of $A^{(j)}$ associated with eigenmode $\phi$ of $A$ by
$\mu^{(1)}\mu^{(2)}\cdots\mu^{(j)}$.
The higher-order in $\delta$ terms of \eqref{eqn:lem-by-mode} are given in a 
form that will be used in the next technical lemma. The $\cP^i()$ notation is introduced to state the relevant result without keeping track of the specific factors in each product.

\begin{proof}
The proof relies on two repeated calculations. First, for any $\alpha,\beta$
\begin{align}\label{eqn:ldc-calc1}
A_\beta \begin{pmatrix} \alpha \phi \\ 0 \end{pmatrix}= 
  \begin{pmatrix} A & -\beta I \\ I & 0 \end{pmatrix}
       \begin{pmatrix} \alpha \phi \\ 0 \end{pmatrix} = 
  \alpha \begin{pmatrix} \lambda \phi \\ \phi \end{pmatrix} = 
  \alpha \begin{pmatrix} \mu^{(0)} \phi \\ \phi \end{pmatrix},
\end{align}
where $\mu^{(0)}: = \lambda$.
Second, 
\begin{align}\label{eqn:ldc-calc2}
A^{(k)} \begin{pmatrix} \mu^{(j)} \phi \\ \phi \end{pmatrix} &= 
A^{(k)} \begin{pmatrix} \mu^{(k)} \phi \\ \phi \end{pmatrix} + 
A^{(k)} \begin{pmatrix} \delta_{jk} \phi \\ 0 \end{pmatrix} = 
\mu^{(k)} \begin{pmatrix} \mu^{(k)} \phi \\ \phi \end{pmatrix} + 
A^{(k)} \begin{pmatrix} \delta_{jk} \phi \\ 0 \end{pmatrix} 
\nonumber \\ & =
\mu^{(k)} \begin{pmatrix} \mu^{(k)} \phi \\ \phi \end{pmatrix} + 
\delta_{jk}  \begin{pmatrix} \mu^{(0)} \phi \\ \phi \end{pmatrix},
\end{align}
where the last term in \eqref{eqn:ldc-calc2} is the result of \eqref{eqn:ldc-calc1}.

Starting with  
\eqref{eqn:ldc-calc1}, and 
proceeding to apply \eqref{eqn:ldc-calc2} we have
\begin{align}\label{eqn:ldc-base1}
A^{(0)} & \begin{pmatrix} \phi \\ 0 \end{pmatrix} = 
\begin{pmatrix} \mu^{(0)}\phi \\ \phi \end{pmatrix},
\\
A^{(1)} A^{(0)}& \begin{pmatrix} \phi \\ 0 \end{pmatrix} = 
\mu^{(1)} \begin{pmatrix} \mu^{(1)} \phi \\ \phi \end{pmatrix} + 
\delta_{01}  \begin{pmatrix} \mu^{(0)} \phi \\ \phi \end{pmatrix},
\\
A^{(2)}A^{(1)} A^{(0)}& \begin{pmatrix} \phi \\ 0 \end{pmatrix} = 
\mu^{(2)}\mu^{(1)} \begin{pmatrix} \mu^{(2)} \phi \\ \phi \end{pmatrix} + 
\delta_{12}\mu^{(1)}  \begin{pmatrix} \mu^{(0)} \phi \\ \phi \end{pmatrix} +
\delta_{01}\mu^{(2)} \begin{pmatrix} \mu^{(2)} \phi \\ \phi \end{pmatrix} + 
\delta_{02}\delta_{01}  \begin{pmatrix} \mu^{(0)} \phi \\ \phi \end{pmatrix}
\nonumber \\ & = 
\left( \mu^{(2)}\mu^{(1)} + \delta_{01}\mu^{(2)} \right)
\begin{pmatrix} \mu^{(2)} \phi \\ \phi \end{pmatrix} + 
\left( \delta_{12}\mu^{(1)} + \cP^2(\delta) \right) 
\begin{pmatrix} \mu^{(0)} \phi \\ \phi \end{pmatrix}.
\end{align}
One more iteration reveals the form of the higher order terms.
\begin{align}\label{eqn:ldc-base2}
A^{(3)}A^{(2)}A^{(1)} A^{(0)} \begin{pmatrix} \phi \\ 0 \end{pmatrix} &= 
\left( \mu^{(2)}\mu^{(1)} + \delta_{01}\mu^{(2)} \right)
\left( \mu^{(3)} \begin{pmatrix} \mu^{(3)} \phi \\ \phi \end{pmatrix} + 
\delta_{23}  \begin{pmatrix} \mu^{(0)} \phi \\ \phi \end{pmatrix}\right)
\nonumber \\ &+
\left( \delta_{12}\mu^{(1)} + \cP^2(\delta) \right) 
\left( \mu^{(3)} \begin{pmatrix} \mu^{(3)} \phi \\ \phi \end{pmatrix} + 
\delta_{03}  \begin{pmatrix} \mu^{(0)} \phi \\ \phi \end{pmatrix} \right)
\nonumber \\ & = 
\left( \mu^{(3)} \mu^{(2)} \mu^{(1)} + \delta_{01}\mu^{(3)} \mu^{(2)} 
       + \delta_{12}\mu^{(3)}\mu^{(1)} + 
\cP^2(\delta) \cP(\mu)  \right)
\begin{pmatrix} \mu^{(3)} \phi \\ \phi \end{pmatrix}
\nonumber \\ & + 
\left( \delta_{23} \mu^{(2)} \mu^{(1)} + 
\cP^2(\delta)\cP( \mu) + \cP^3(\delta)\right)
\begin{pmatrix} \mu^{(0)} \phi \\ \phi \end{pmatrix}.
\end{align}

Now we may proceed inductively.  Suppose 
\begin{align}\label{eqn:ldc-ind0}
\Phi^{(j)} := A^{(j)}\cdots A^{(0)} \begin{pmatrix} \phi \\ 0 \end{pmatrix}
= \left( \prod_{i = 1}^j \mu^{(i)}  +
   \sum_{k=1}^{j-1} \delta_{k-1,k}\prod_{i = 1, i \ne k}^j \mu^{(i)} 
   + \sum_{i = 2}^{j-1} 
\cP^i(\delta) \cP^{j-i}(\mu)\right)
\begin{pmatrix} \mu^{(j)} \phi \\ \phi \end{pmatrix}
\nonumber \\
+ \left( \delta_{j-1,j} \prod_{i=1}^{j-1}\mu^{(i)} 
+ \sum_{i = 2}^j \cP^i(\delta)\cP^{j-i}(\mu) \right)
\begin{pmatrix} \mu^{(0)} \phi \\ \phi \end{pmatrix}.
\end{align}
We will show 
\begin{align}\label{eqn:ldc-ind1}
\Phi^{(j+1)} 
& = \left( \prod_{i = 1}^{j+1} \mu^{(i)}  +
  \sum_{k=1}^{j} \delta_{k-1,k}\prod_{i = 1, i \ne k}^{j+1} \mu^{(i)} 
 +\sum_{i = 2}^j \cP^i(\delta) \cP^{j+1-i}(\mu)
  \right)
\begin{pmatrix} \mu^{(j+1)} \phi \\ \phi \end{pmatrix}
\nonumber \\ & 
+ \left( \delta_{j,j+1} \prod_{i=1}^{j}\mu^{(i)} + \sum_{i = 2}^{j+1}
\cP^i(\delta) \cP^{j+1-i}\mu) \right)
\begin{pmatrix} \mu^{(0)} \phi \\ \phi \end{pmatrix}.
\end{align}

The base step of the induction is satisfied by \eqref{eqn:ldc-base2}. For the inductive
step, applying \eqref{eqn:ldc-calc2} to $A^{(j+1)}\Phi^{(j)}$ yields
\begin{align*}
A^{(j+1)}\Phi^{(j)} &= 
\mu^{(j+1)}\left\{
  \left( \prod_{i = 1}^j \mu^{(i)}  +
   \sum_{k=1}^{j-1} \delta_{k-1,k}\prod_{i = 1, i \ne k}^j \mu^{(i)} 
  + \sum_{i = 2}^{j-1} \cP^i(\delta) \cP^{j-i}(\mu) \right) \right. 
  \nonumber \\ &+ \left.
  \left( \delta_{j-1,j} \prod_{i=1}^{j-1}\mu^{(i)} 
 + \sum_{i = 2}^j \cP^i(\delta) \cP^{j-i}(\mu)  
\right) \right\}
\begin{pmatrix} \mu^{(j+1)} \phi \\ \phi \end{pmatrix} 
\nonumber \\
 &+ \left\{ \delta_{j,j+1}
  \left( \prod_{i = 1}^j \mu^{(i)}  +
   \sum_{k=1}^{j-1} \delta_{k-1,k}\prod_{i = 1, i \ne k}^j \mu^{(i)} 
 + \sum_{i = 2}^{j-1} \cP^i(\delta)\cP^{j-i}(\mu)
   \right)
\right. \nonumber \\ & \left.
  + \delta_{0,j+1}
  \left( \delta_{j-1,j} \prod_{i=1}^{j-1}\mu^{(i)} 
+ \sum_{i = 2}^j \cP^i(\delta) \cP^{j-i}(\mu) 
\right) \right\}
\begin{pmatrix} \mu^{(0)} \phi \\ \phi \end{pmatrix},
\end{align*}
which after multiplying though and combining the $\cP( \cdot )$ terms of like order 
agrees with \eqref{eqn:ldc-ind1}.
\end{proof}
This establishes the result \eqref{eqn:lem-by-mode}.
%% -- --------------------------------------------------------------------
The next step in the argument is to generalize the first component of the initial vector
used in lemma \ref{lem:dm-by-mode} from a single eigenmode of $A$ to  
a linear combination of eigenmodes of $A$, to arrive an an estimate analogous
to \eqref{eqn:diagpow}.   

%% -- --------------------------------------------------------------------
\begin{lemma}\label{lem:dyprod}
Let $A$ satisfy assumption \ref{assume:l1diag}.
Let $\delta_{l,i,k} = \mu_l^{(i)} - \mu_l^{(k)}$. 
As in lemma \ref{lem:dm-by-mode}, 
define $\cP^{i}(\mu_l)$ to be a product of $i$ terms $\mu_l^{(k)}$, 
where $1 \le k \le j$, 
and $\cP^i(\delta_l)$ to be a product of $i$ terms $\delta_{l,k,p}$, where 
$0 \le k,p \le j$.
Let $u_0 = \sum_{l = 1}^n a_l \phi_l$, a linear combination of the eigenvectors of $A$.
Then it holds that the product 
$A^{(j)}A^{(j-1)}\cdots A^{(0)}\begin{pmatrix}u_0\\0\end{pmatrix}$ satisfies
\begin{align}\label{eqn:lem-dyprod}
A^{(j)}A^{(j-1)}\cdots A^{(0)} \begin{pmatrix} u_0 \\ 0 \end{pmatrix}
& = a_1 \left( \prod_{i = 1}^j \mu_1^{(i)} \right) \left\{
\psi_1^{(j)} + \sum_{l = 2}^n \f{a_l}{a_1} 
\left( \prod_{i = 1}^j\f{\mu_l^{(i)}}{\mu_1^{(i)}} \right)
\psi_l^{(j)} \right\}
\nonumber \\ & +
\sum_{i = 1}^{j-1} a_1 \cP^{j-i}(\mu_1) 
\left\{ \psi_1^{(j)} \cP^i(\delta_1) + \sum_{l=2}^n \f{a_l}{a_1} 
 \cP^i(\delta_l)  \cP^{j-i}\left(\f{\mu_l}{\mu_1} \right) \psi_l^{(j)} 
\right\}
\nonumber \\ & +
\sum_{i = 1}^j  a_1 \cP^{j-i}(\mu_1) 
\left\{
\psi_1^{(0)} \cP^i(\delta_1) + \sum_{l = 2}^n \f{a_l}{a_1} 
\cP^i(\delta_l) \cP^{j-i} \left(\f{\mu_l}{\mu_1} \right)  \psi_l^{(0)}  
\right\}.
\end{align}

Supposing additionally that $\mu_1^{(i)} > \delta_{l,k,p}$ for any $i,k,p = 1, \ldots, j$,
and $l \ge 2$,
then as $j$ increases, the product $A^{(j)}A^{(j-1)}\cdots A^{(0)}u_0$ aligns to a 
linear combination of $\psi_1^{(j)}$ and $\psi_1^{(0)}.$

\end{lemma}
The proof shows additional detail on the $\bigo(\delta)$ terms, as revealed in lemma
\ref{lem:dm-by-mode}. This lemma shows that the product of the sequence of matrices
$A^{(j)}\ldots A^{(0)}$ applied to a vector with a general first component, 
$u_0 \in \R^n$ and null second component $0 \in \R^n$ 
aligns with a vector whose first component is
the dominant eigenvector $\phi_1$ of $A$.
It will be shown in theorem \ref{thm:dymcon} that the convergence is similar to the 
power method with $(\lambda_l/\lambda_1)^j$ as in \eqref{eqn:diagpow}
replaced by the product
$(\mu_l^{(1)}\cdots \mu_l^{(j)})/(\mu_1^{(1)}\cdots\mu_1^{(j)})$.
The appreciable difference in the convergence is from the contribution of the 
$\delta$-scaled terms which are in 
the directions of the eigenvectors $\psi_l^{(j)}$ and $\psi_l^{(0)}$, with 
$l = 1, \ldots, n$. 
As we will see in theorem \ref{thm:dymcon} and remark \ref{rem:resicoef}, 
these terms will not interfere with 
convergence or the asymptotically expected rate, due to the stability of the
parameters $\beta_i$, as shown in lemma \ref{lem:rstab}. 

\begin{proof}
First by applying linearity and \eqref{eqn:lem-by-mode} we have
\begin{align}\label{eqn:dyprod01}
A^{(j)}\cdots A^{(0)} \begin{pmatrix} u_0 \\ 0 \end{pmatrix} & =
\sum_{l = 1}^n a_l 
A^{(j)}\cdots A^{(0)} \begin{pmatrix} \phi_l \\ 0 \end{pmatrix}
\nonumber \\ &
= \sum _{l = 1}^n a_l
 \left( \prod_{i = 1}^j \mu_l^{(i)}  +
   \sum_{k=1}^{j-1} \delta_{l,k-1,k}\prod_{i = 1, i \ne k}^j \mu_l^{(i)} 
   + \sum_{i = 2}^{j-1} \cP^i(\delta_l) \cP^{j-i}(\mu_l)\right)
\begin{pmatrix} \mu_l^{(j)} \phi_l \\ \phi_l \end{pmatrix}
\nonumber \\ &
+ \sum _{l = 1}^n a_l
\left( \delta_{l,j-1,j} \prod_{i=1}^{j-1}\mu_l^{(i)} + \sum_{i = 2}^j 
\cP^i(\delta_l) \cP^{j-i}(\mu_l) \right)
\begin{pmatrix} \mu_l^{(0)} \phi_l \\ \phi_l \end{pmatrix}.
\end{align}
Now we will examine each term of \eqref{eqn:dyprod01}.

We rewrite the first term in the right hand side of \eqref{eqn:dyprod01} as
\begin{align}\label{eqn:dyprod02}
\sum_{l = 1}^n a_l \left(\prod_{i = 1}^j \mu_l^{(i)} \right) \psi_l^{(j)}
= a_1 \left( \prod_{i = 1}^j \mu_1^{(i)} \right) \left\{
\psi_1^{(j)} + \sum_{l = 2}^n \f{a_l}{a_1} 
\left( \prod_{i = 1}^j\f{\mu_l^{(i)}}{\mu_1^{(i)}} \right)
\psi_l^{(j)} \right\}.
\end{align}
This is similar to \eqref{eqn:diagpow} 
and displays convergence to $\psi_1$ so long as the other terms do not interfere.
The second term in the right hand side of \eqref{eqn:dyprod01} can be written as
\begin{align}\label{eqn:dyprod03}
&\sum_{l = 1}^n a_l \sum_{k = 1}^{j-1} \delta_{l,k-1,k} 
\left( \prod_{i = 1, i \ne k}^j \mu_l^{(i)} \right)
\psi_l^{(j)} 
\nonumber \\  = &
\sum_{k = 1}^{j-1} a_1 \left(\prod_{i = 1,i \ne k}^j \mu_1^i\right)
\left\{ \delta_{1,k-1,k}\psi_1^{(j)} + 
\sum_{l = 2}^n \f{a_l}{a_1} \delta_{l, k-1,k} \left( \prod_{i = 1, i \ne k}^j 
\f{ \mu_l^{(i)} }{ \mu_{1}^{(i)} } \right) \psi_l^{(j)}
\right\},
\end{align}
which is an $\bigo(\delta)$ term 
where the factors of $\mu^{(i)}_l/\mu^{(i)}_1$ multiplying the
subdominant eigenmodes are one power lower
than in the dominant term \eqref{eqn:dyprod02}.  
The higher order terms multiplying the eigenvectors of $A^{(j)}$ are
\begin{align}\label{eqn:dyprod04}
\sum_{l=1}^n a_l \sum_{i = 2}^{j-1}\cP^i(\delta_l) \cP^{j-i}(\mu_l) \psi_l^{(j)}
=
\sum_{i = 2}^{j-1} a_1 \cP^{j-i}(\mu_1) 
\left\{ \psi_1^{(j)} \cP^i(\delta_1) + \sum_{l=2}^n \f{a_l}{a_1} 
\cP^i( \delta_l)\cP^{j-i}\left(\f{\mu_l}{\mu_1} \right)\psi_l^{(j)} 
\right\}. 
\end{align}

Next, we look at the terms of \eqref{eqn:dyprod01} multiplying the eigenvectors 
$\psi_0^{(0)}$ of $A^{(0)}$. The lowest order term is $\bigo(\delta)$ and is given by
\begin{align}\label{eqn:dyprod05}
\sum_{l = 1}^n a_l \delta_{l,j-1,j} \left( \prod_{i = 1}^{j-1} \mu_l^{(i)} \right)
= a_1  \left( \prod_{i = 1}^{j-1} \mu_1^{(i)} \right) 
\left\{
\delta_{1,j-1,j} \psi_1^{(0)} + \f{a_l}{a_1} \delta_{l,j-1,j} 
\left( \prod_{i = 1}^{j-1} \f{\mu_l^{(i)}}{\mu_1^{(i)}} \right) \psi_l^{(0)}
\right\}.
\end{align}
Last we have the higher order terms
\begin{align}\label{eqn:dyprod06}
\sum_{l = 1}^n a_l \left( \sum_{i = 2}^j \cP^i(\delta_l) \cP^{j-i}(\mu_l) \right) \psi_l^{(0)}
= \sum_{i = 2}^j  a_1 \cP^{j-i}(\mu_1) 
\left\{
\psi_1^{(0)} \cP^i(\delta_1) + \sum_{l = 1}^n \f{a_l}{a_1} 
\cP^i(\delta_l) \cP^{j-i}\left(\f{\mu_l}{\mu_1} \right) \psi_l^{(0)}  
\right\}.
\end{align}

Sweeping the results of the more detailed \eqref{eqn:dyprod03} into 
\eqref{eqn:dyprod04}, and likewise
\eqref{eqn:dyprod05} into \eqref{eqn:dyprod06} yields the result \eqref{eqn:lem-dyprod}.
Finally, the alignment of the product \eqref{eqn:lem-dyprod} to a combination
of $\psi_1^{(j)}$ and $\psi_1^{(0)}$ follows from noting each ratio
$(\mu^{(i)}_l/\mu^{(i)}_1) < 1$ and applying the hypothesis   
$\mu_1^{(i)} > \delta_{l,k,p}$ for any $i,k,p = 1, \ldots, j$,
and $l \ge 2$.
\end{proof}
%% -- --------------------------------------------------------------------
%% -- --------------------------------------------------------------------
The next lemma shows that if $\rho_k$ is an $\eps$ perturbation of $\rho = |\mu_2/\mu_1|$,
then $r_{k+1}$ is an $\widehat \eps$ perturbation of $r$, where $\widehat \eps < 2\eps$
for $\rho \in (0,1)$, and $\widehat \eps \goto 0$ as $\rho \goto 1$. This means the
smaller the spectral gap in $A$, the more stable the dynamic momentum method becomes. 
\begin{lemma}\label{lem:rstab} 
Let $\rho \in (0,1)$ and consider $\eps$ small enough so 
that $(2 \rho \eps + \eps^2)/(1+\rho^2) < 1$.
Let $\rho_k = \rho + \eps$ and define 
$r_{k+1} = 2 \rho_k/(1 + \rho_k^2)$, as in algorithm \ref{alg:dymo}. 
Then
\begin{align}\label{eqn:lemrstab}
r_{k+1} = r + \widehat \eps + \bigo(\eps^2) ~\text{ with }~
\widehat \eps = \eps \f{2 (1-\rho^2)}{(1+\rho^2)^2}.
\end{align}
\end{lemma}
The condition $(2 \rho \eps + \eps^2)/(1+\rho^2) < 1$ is satisfied for $\rho \in (0,1)$ by
$\eps < 0.71$.

\begin{proof}
For $r = |\lambda_2/\lambda_1|$ the asymptotic convergence rate of iteration 
\eqref{eqn:mom1} is
$\rho = r(1 + \sqrt{1-r^2})^{-1}$, as given by \eqref{eqn:cor-momrate}, when 
$\beta = \lambda_2^2/4$. Inverting this expression for $r$ in terms of $\rho$ yields
\begin{align}\label{eqn:rofrho}
r = \f{2\rho}{1 + \rho^2}.
\end{align}
Suppose the detected convergence rate $\rho_k = d_{k+1}/d_k$ is an $\epsilon$ 
perturbation of $\rho$, meaning $\rho_k = \rho + \eps$.
Expanding $r_{k+1} = 2 \rho_k/(1+\rho_k^2)$ in $\eps$ yields
\begin{align}\label{eqn:beta01}
r_{k+1} & = \f{2 \rho}{1 + (\rho + \eps)^2} + \f {2\eps}{1 + (\rho + \eps)^2}
\nonumber \\
& = \f{2\rho}{ 1 + \rho^2} \left( \f{1}{1 + \f{2 \rho \eps + \eps^2}{1 + \rho^2}}\right)
+  \f{2\eps}{ 1 + \rho^2} \left( \f{1}{1 + \f{2 \rho \eps + \eps^2}{1 + \rho^2}}\right)
\nonumber \\ 
& = \f{2\rho}{ 1 + \rho^2} \left( 1 - \f{2 \rho \eps + \eps^2}{1 + \rho^2} 
+ \bigo(\eps^2)\right)
+ \f{2 \eps}{1 + \rho^2} + \bigo(\eps^2). 
\end{align}
Applying \eqref{eqn:rofrho} to \eqref{eqn:beta01} yields
\begin{align}\label{eqn:beta02}
r_{k+1} = r\left(1 + \eps \left( \f 1 \rho - r\right) \right) + \bigo(\eps^2) 
= r +  \widehat \eps + \bigo(\eps^2),
~ \text{ where }~
\widehat \eps = r\eps \left( \f 1 \rho - r\right).
\end{align}
Applying \eqref{eqn:rofrho} to \eqref{eqn:beta02} yields
the result \eqref{eqn:lemrstab},
by which  $\widehat \eps < 2\eps$ for $\rho \in (0,1)$, 
and $\widehat \eps < \eps$ for $\rho \in (0.486,1)$, or $r \in (0.786,1)$.
Moreover as $r$ getting closer to unity, the approximation becomes more stable,
with $\widehat \eps < 0.161 \cdot \eps$ for $r \in (0.99,1)$ and 
     $\widehat \eps < 0.0468 \cdot \eps$ for $r \in (0.999,1)$. 
\end{proof}
The stability of $\beta_k = (r_k \nu_k)^2/4 $ in algorithm \ref{alg:dymo}
is inherited directly from the
stability of $r_k$, once $\nu_k$ sufficiently converges to $\lambda_1$.
%% -- --------------------------------------------------------------------
%% -- --------------------------------------------------------------------
\begin{remark}\label{rem:betastab}
Another way to view how close $\beta_k$ is to $\beta = \beta_{opt} = \lambda_2^2/4$ 
with respect to 
$r_k$ and $\rho_k$ viewed as perturbations of $r$ and $\rho$ is to consider 
$\rho_k$ written as
\[
\rho_k = \f{r \sqrt{1 + \eps/r^2}}{ 1 + \sqrt{1 - r^2(1 + \eps/r^2)}},
\]
for some $\eps$ with  $-r^2 < \eps < 1-r^2$.
Applying $r_{k+1} = 2 \rho_k/(1+\rho_k^2)$ we then have
$
r_{k+1} = r \sqrt{1 + \eps/r^2},
$
by which  $\beta_{k+1} = r^2(1 + \eps/r^2) \nu_{k+1}^2/4$. For $\nu_{k+1} \approx \lambda_1$
this yields
\[
\beta_{k+1} \approx \f{\lambda_2^2}{4} + \eps \f{\lambda_1^2}{4},
\]
which shows how perturbations $r_k$ with respect to $r$ result in perturbations to 
$\beta_k$ with respect to $\beta$.
\end{remark}
%% -- --------------------------------------------------------------------
%% -- --------------------------------------------------------------------
Now we can summarize the results of this section in a convergence theorem.
\begin{theorem}\label{thm:dymcon}
Let $A$ satisfy assumption \ref{assume:l1diag}. 
Let $\delta_{l,i,k} = \mu_l^{(i)} - \mu_l^{(k)}$, and 
let $u_0 = \sum_{l = 1}^n a_j \phi_j$, 
a linear combination of the eigenvectors of $A$.
If $A$ is symmetric then $\beta_k < \lambda_1^2/4$ for all $k$.
Then \eqref{eqn:lem-dyprod} holds and algorithm \ref{alg:dymo} converges to 
the dominant eigenpair.
\end{theorem}
Here we proceed by assuming generically that none of the $\beta_k$ take a value of 
exactly equal to 
$\lambda^2/4$ for any eigenvalue $\lambda$ of $A$. This is a reasonable assumption
due both to floating point arithmetic, and that as shown in lemma \ref{lem:rstab}, 
the $\beta_k$ only converge to $\lambda_2^2/4$ as $r \goto 1$, and we are always
in the circumstance that $r < 1$.
\begin{proof}
By the definitions of $\rho_k$ and $r_k$, we have $r_k^2 \le 1$.  Since $\nu_k$ 
is the Rayleigh quotient with approximate eigenvector $x_{k+1}$ and symmetric $A$, 
it follows that $\beta_k < \lambda_1^2/4$.
We will start by developing bounds on the $\mu_l^{(i)}$ and the 
$\delta_{l,i,k}$, and in the process will verify the final hypothesis of
lemma \ref{lem:dyprod} by verifying $|\delta_{l,i,k}| \le \max\{|\mu_l^{(i)}|,|\mu_l^{(k)}|\}$.
We will also see that $\delta_{l,i,k} \goto 0$ as $\beta_i - \beta_k \goto 0$ for each
$\lambda_l$.
Consider $\lambda = \lambda_l \ne 0$. 
There are three cases we need to consider. Without loss of generality, suppose 
$\beta_i \ge \beta_k$. 
\begin{enumerate}
\item[(i)] If $\lambda^2/4 \ge \beta_i$ then $\mu_l^{(i)}$ and 
$\mu_l^{(k)}$ are given by \eqref{eqn:mucase1}, and
$|\mu_l^{(i)}| \in [|\lambda|/2,|\lambda|]$ for $l = 2, \ldots, n$. 
Since we have $\beta < \lambda_1^2/4$ we have $|\mu_1^{(i)}| \in (|\lambda_1|/2, |\lambda_1|]$. 
To bound $\delta_{l,i,k}$ we have
\begin{align}\label{eqn:ct001}
|\mu_l^{(i)} - \mu_{l}^{(k)}| = \f{|\lambda|}{2}
\left|\sqrt{1 - 4\beta_i/\lambda^2} - \sqrt{1 - 4\beta_k/\lambda^2} \right| < \f{|\lambda|}{2}
\le |\mu_l^{(i)}| < |\mu_l^{(k)}|.
\end{align}
It is clear from continuity and \eqref{eqn:ct001} that $\delta_{l,i,k}\goto 0$ as 
$\beta_i-\beta_k \goto 0$. 
We can also expand to first order to see how, yielding
$|\delta_{l,i,k}| = |(\beta_i-\beta_k) + \ldots|$.

\item[(ii)] If $\lambda^2/4 \le \beta_i, \beta_k < \lambda_1^2/4$ then 
$\mu_l^{(i)}$ and $\mu_l^{(k)}$ are given by \eqref{eqn:mucase3}, and
we have $|\mu_l^{(i)}| = \sqrt{\beta_i}$.
In this case $\delta_{l,i,k}$ satisfies
\begin{align}\label{eqn:ct002}
 |\mu_l^{(i)} - \mu_{l}^{(k)}| = 
\left| 
\sqrt{\beta_i}\sqrt{1 - \lambda^2/(4\beta_i)} -       
\sqrt{\beta_k}\sqrt{1 - \lambda^2/(4\beta_k)} 
\right| 
\le \sqrt{\beta_i} = |\mu_l^{(i)}|.
\end{align}
From continuity and \eqref{eqn:ct002} it is clear that $\delta_{l,i,k}\goto 0$ as
$\beta_k - \beta_i \goto 0$. 
Expanding \eqref{eqn:ct002} to first order to see how, yields
$|\mu_l^{(i)} - \mu_{l}^{(k)}| =|(\sqrt{\beta_i} - \sqrt{\beta_k})(1 - \lambda^2/8)+\ldots|$.

\item[(iii)] If $\beta_k \le \lambda^2/4 \le \beta_i$, then we have 
\[
\mu_l^{(i)} - \mu_l^{(k)} = 
\f{\lambda}{2} + \f 1 2 \sqrt{\lambda^2 - 4 \beta_i} - 
\left( \f{\lambda}{2} + \f 1 2 \sqrt{\lambda^2 - 4 \beta_k}\right) =
\f{i}{2} \sqrt{4 \beta_i - \lambda^2} - \f 1 2\sqrt{\lambda^2 - 4 \beta_k},
\]
by which $|\mu_l^{(i)} - \mu_l^{(k)}| = \sqrt{\beta_i - \beta_k} \le \sqrt{\beta_i}
= |\mu_l^{(i)}|$.
\end{enumerate}
Combining with the above results, we have $|\delta_{l,k,p}| \le 
\max\{|\mu_l^{(k)}|, |\mu_l^{(p)}|\}$, for any $l,k,p = 1, \ldots, j$.

We now have by lemma \ref{lem:dyprod} that as $j$ increases, 
the product $A^{(j)}A^{(j-1)}\cdots A^{(0)}\begin{pmatrix}u_0 \\ 0 \end{pmatrix}$
aligns with a linear combination of $\psi_1^{(j)}$ and $\psi_1^{(0)}$. 
As in subsection \ref{subsec:static}, we now analyze the convergence of algorithm~
\ref{alg:dymo} by the convergence of 
\begin{align}\label{eqn:ct003}
\begin{pmatrix}x_{j+1}\\y_{j+1} \end{pmatrix} = 
\f 1 {h_{j+1}}\begin{pmatrix}u_{j+1}\\z_{j+1} \end{pmatrix} = 
\f 1 {h_{j+1}}
\left( \prod_{i = 0}^{j} h_i^{-1} \right) A^{(j)}A^{(j-1)}\cdots A^{(0)}
\begin{pmatrix}u_0 \\ 0 \end{pmatrix},
\end{align}
where $h_j = \nr{u_j}$. 
Applying \eqref{eqn:lem-dyprod} to \eqref{eqn:ct003}, we have
\begin{align}\label{eqn:ct004}
\begin{pmatrix}u_{j+1}\\z_{j+1} \end{pmatrix} 
& = \f{a_1}{h_0} \left(\prod_{i = 1}^j \f{\mu_1^{(i)}}{h_i} \right) \left\{
\psi_1^{(j)} + \sum_{l = 2}^n \f{a_l}{a_1} 
\left( \prod_{i = 1}^j\f{\mu_l^{(i)}}{\mu_1^{(i)}} \right)
\psi_l^{(j)} \right\}
\nonumber \\ & +
\left( \prod_{i = 0}^j \f {1}{h_i}\right)
\sum_{i = 1}^{j-1} a_1 \cP^{j-i}(\mu_1) 
\left\{ \psi_1^{(j)} \cP^i(\delta_1) + \sum_{l=2}^n \f{a_l}{a_1} 
\cP^i(\delta_l) \cP^{j-i} \left( \f{\mu_l}{\mu_1} \right) \psi_l^{(j)} 
\right\}
\nonumber \\ & +
\left( \prod_{i = 0}^j \f {1}{h_i} \right)
\sum_{i = 1}^j  a_1 \cP^{j-i} (\mu_1) 
\left\{
\psi_1^{(0)} \cP^i(\delta_1) + \sum_{l = 2}^n \f{a_l}{a_1} 
\cP^i ( \delta_l) \cP^{j-i}\left(\f{\mu_l}{\mu_1} \right) \psi_l^{(0)}  
\right\}.
\end{align}
Distributing through the normalization factors in \eqref{eqn:ct004} yields
\begin{align}\label{eqn:ct005}
\begin{pmatrix}u_{j+1}\\z_{j+1} \end{pmatrix}
& =  
 \f{a_1}{h_0} \left(\prod_{i = 1}^j \f{\mu_1^{(i)}}{h_i} \right) \left\{
\psi_1^{(j)} + \sum_{l = 2}^n \f{a_l}{a_1} 
\left( \prod_{i = 1}^j\f{\mu_l^{(i)}}{\mu_1^{(i)}} \right)
\psi_l^{(j)} \right\}
\nonumber \\ & +
\sum_{i = 1}^{j-1} \f{a_1}{h_0} 
\left( \f{\cP^{j-i}(\mu_1)}
{\prod_{k = 1}^{j-i} h_k} \right)
\left\{ \psi_1^{(j)} \cP^i(\delta_1) + \sum_{l=2}^n \f{a_l}{a_1} 
\left( \f{\cP^i( \delta_l)}{\prod_{k = j-i+1}^{j} h_k}\right)
\cP^{j-i}\left(\f{\mu_l}{\mu_1} \right)\psi_l^{(j)} 
\right\}
\nonumber \\ & +
\sum_{i = 1}^j  \f{a_1}{h_0} \left( \f{\cP^{j-i} (\mu_1)} {\prod_{k = 1}^{j-i} h_k} \right)
\left\{
\psi_1^{(0)} \cP^i(\delta_1) + \sum_{l = 2}^n \f{a_l}{a_1} 
\left( \f{\cP^i( \delta_l)}{\prod_{k = j-i+1}^{j} h_k}\right)
\cP^{j-i} \left(\f{\mu_l}{\mu_1} \right) \psi_l^{(0)}  
\right\}.
\end{align}

By the arguments above, 
$\begin{pmatrix}u_{j+1}\\z_{j+1} \end{pmatrix}$
aligns with a linear combination of $\psi_1^{(j)}$ and $\psi_1^{(0)}$,
both of which have first components in the direction of $\phi_1$.
This further shows that the Rayleigh quotient $\nu_k=  (Ax_k,x_k) \goto \lambda_1$, by 
which the residual \eqref{eqn:resi} converges to zero.
\end{proof}
%% -- --------------------------------------------------------------------
We conclude this section with a heuristic discussion of the coefficients of each eigenmode
that appear in the residual, as per remark \ref{rem:resi}.
%% -- --------------------------------------------------------------------
\begin{remark}\label{rem:resicoef}
By theorem \ref{thm:dymcon}, the Rayleigh quotient $\nu_k$ converges to $\lambda_1$.
As in remark \ref{rem:resi}, we consider the asymptotic regime where 
$\nu_k \approx \lambda_1$, so that the ratio between consecutive residuals 
$\rho_k$ is well approximated by
\begin{align}\label{eqn:ct005a}
\rho_j \approx 
\f{\nr{\sum_{l = 2}^n (\lambda_l - \nu_{j+1}) \alpha_l^{(j+1)} \phi_l}}
{\nr{\sum_{l = 2}^n (\lambda_l - \nu_j) \alpha_l^{(j)} \phi_l}}.
\end{align}
From \eqref{eqn:ct005a}, and the definition of the eigenvectors of the augmented matrix in 
\eqref{eqn:aug-vect}, the coefficients $\alpha_l^{(j+1)}$, $l \ge 2$  are given by
\begin{align}\label{eqn:ct006}
\alpha_l^{(j+1)} =  
 \f{a_l}{\prod_{i = 0}^{j+1} h_{i}} \left\{ \mu_l^{(j)} \left(\prod_{i = 1}^j \mu_l^{(i)}  \right)
 +
\sum_{i = 1}^{j-1} (\mu_l^{(j)} + 1)
\cP^{j-i}(\mu_l) \cP^i(\delta_l) 
+  \cP^j(\delta_l) \right\}.
\end{align}
We next make the argument that the first term inside the brackets in \eqref{eqn:ct006} 
dominates the others.

From theorem \ref{thm:dymcon}, each $\delta_{l,i,k}$ satisfies 
$|\delta_{l,i,k}| \le \max\{\mu_l^{(i)}, \mu_l^{(k)}\}$.
Referring to the proof of lemma \ref{lem:dm-by-mode}, each of the factors of 
$\delta_{l,i,k}$ have
either the form $\delta_{l,p-1,p}$ or $\delta_{l,0,p}$, where $p$ ranges from $1$ to $j+1$.  
As per the discussion in theorem \ref{thm:dymcon}, 
the terms of the form $\delta_{l,p-1,p}$ go to zero as the 
$\beta_k \goto \beta = \lambda_2^2/4$.   
By lemma \ref{lem:rstab}, considering the detected convergence rate $\rho_k$ as a 
perturbation of
the theoretically optimal rate $\rho$, the computed approximation $r_{k+1}$ to 
$r = |\lambda_2/\lambda_1|$ 
is restricted to a tighter interval about $r$ when $r$ is closer to one.
By this argument, and remark \ref{rem:betastab}, $\beta_{k+1}$ is restricted to a small 
interval
around $\beta$ (smaller as $r$ getting closer to 1). So as $j$ increases, terms with of the form
$\delta_{l,j-1,j}$ become negligible. By these arguments, each of the terms under the sum 
of \eqref{eqn:ct006} should be of equal order or less than the first term, and as $j$ 
increases, additional terms under the sum should be essentially negligible.

By inspecting the proof of lemma \ref{lem:dm-by-mode}, 
the final term in \eqref{eqn:ct006} can be seen 
to be $\delta_{l,0,1}\delta_{l,0,2}\cdots\delta_{l,0,j}$. 
By the same arguments above, this term should also be of equal order or less than the 
first, although $\delta_{l,0,j}$ is not in general expected to become negligible as $j$ 
increases. In conclusion, the coefficients $\alpha_l^{(j)}$
are dominated by the products of the eigenvalues $\mu_l^{(i)}, \ i = 1, \dots, j$.

\end{remark}

% -- --------------------------------------------------------------------
%% -- --------------------------------------------------------------------
\section{Numerical results}\label{sec:numerics}
In this section we include four suites of tests comparing the introduced dynamic
momentum method algorithm \ref{alg:dymo} with the power method algorithm \ref{alg:pow},
and the static momentum method with optimal $\beta = \lambda_2^2/4$ as in algorithm
\ref{alg:hbpow}. We include additional comparisons in the first three test suites
with the delayed momentum power method (DMPOW), \cite[algorithm 1]{RJRH22}. In the
last test suite we include comparisons with algorithm \ref{alg:hbpow} with 
the parameter $\beta$ replaced by small perturbations above and below the
optimal value. 

In our implementation of DMPOW we do not assume any spectral knowledge, and we 
consider 20, 100 and 500 preliminary power iterations with deflation in the preliminary
stage to determine an 
approximation of $\lambda_2$. As each of the preliminary iterations contains 3 
matrix-vector multiplications, that number where it is reported exceeds the number of
total iterations for DMPOW as it includes both stages of the algorithm.  
The other methods tested
each require one matrix-vector multiply per iteration.  
We found we were able to improve the performance of DMPOW by choosing $w_0$, which is the
the initial approximation to the second eigenvector, to be orthogonal to $u_0$ (denoted $q_0$ in \cite{RJRH22}). 
We used this technique in DMPOW for all reported results.

All tests were performed in Matlab R2023b running on an Apple MacBook Air with 24GB of
memory, 8 core CPU with 8 core GPU. Throughout this section, each iteration was 
run to a maximum of 2000 iterations or a residual tolerance of $10^{-12}$.  
We include tests started from the fixed initial iterate 
$u_0 = \begin{pmatrix}1 & 1 & \cdots & 1\end{pmatrix}^T$ so that the results can
be reproduced, as well as tests starting from random initial guesses via
{\tt u0 = (rand(n,1)-0.5);}.  To run algorithm \ref{alg:hbpow} which requires
knowledge of $\lambda_2$, we recovered the first two eigenvalues using
{\tt eigs(A,2)}. We emphasize that we did this for comparison purposes only, and that
our interest is in developing effective methods that
do not require any a priori information of the spectrum.
%% -- --------------------------------------------------------------------
\subsection{Test suite 1}\label{subsec:test1}

Our first test suite consists of three symmetric positive definite (SPD) benchmark problems.
All three matrices have similar values of $r \approx 0.999$. This first is a diagonal
matrix included for its transparency. 
The second matrix {\tt Kuu} is used as a benchmark in \cite{PoSc21}.
The third, {\tt Muu} features $\lambda_2 = \lambda_1$, so we demonstrate replacing
$\lambda_2$ with $\lambda_3$ in algorithm \eqref{alg:hbpow}. Our dynamic 
algorithm \ref{alg:dymo} works as expected without modification.
\begin{description}
\item[Matrix 1: ] $A = \diag(1000:-1:1)$. This matrix is a standard benchmark with 
$r = 0.999$.
\item[Matrix 2: ] $A = ${\tt Kuu} from \cite{DH11}, with $n = 7102$. 
This matrix has leading eigenvalues
$\lambda_1 = 54.0821$ and $\lambda_2 = 53.9817$, with $r = 0.9981$.
\item[Matrix 3: ] $A = ${\tt Muu} from \cite{DH11}, with $n = 7102$. 
This matrix has leading 
eigenvalues 
$\lambda_1 = 10^{-3} \times 0.8399$, $\lambda_2 = 10^{-3} \times 0.8398$ and 
$\lambda_3 = 10^{-3} \times 0.8391$.
Using {\tt eigs}, $\lambda_1$ and $\lambda_2$ 
agreed to $10^{-14}$, and algorithm \ref{alg:hbpow} did not converge using 
$\beta = \lambda_2^2/4$. The results shown use $\beta = \lambda_3^2/4$, as $\lambda_3$
is the second largest eigenvalue for this matrix. Taking in this case
$r = \lambda_3/\lambda_1$ yields $r = 0.9992$.
\end{description}

\begin{figure}
\includegraphics[trim = 0pt 0pt 0pt 0pt,clip = true, width = 0.32\textwidth]
{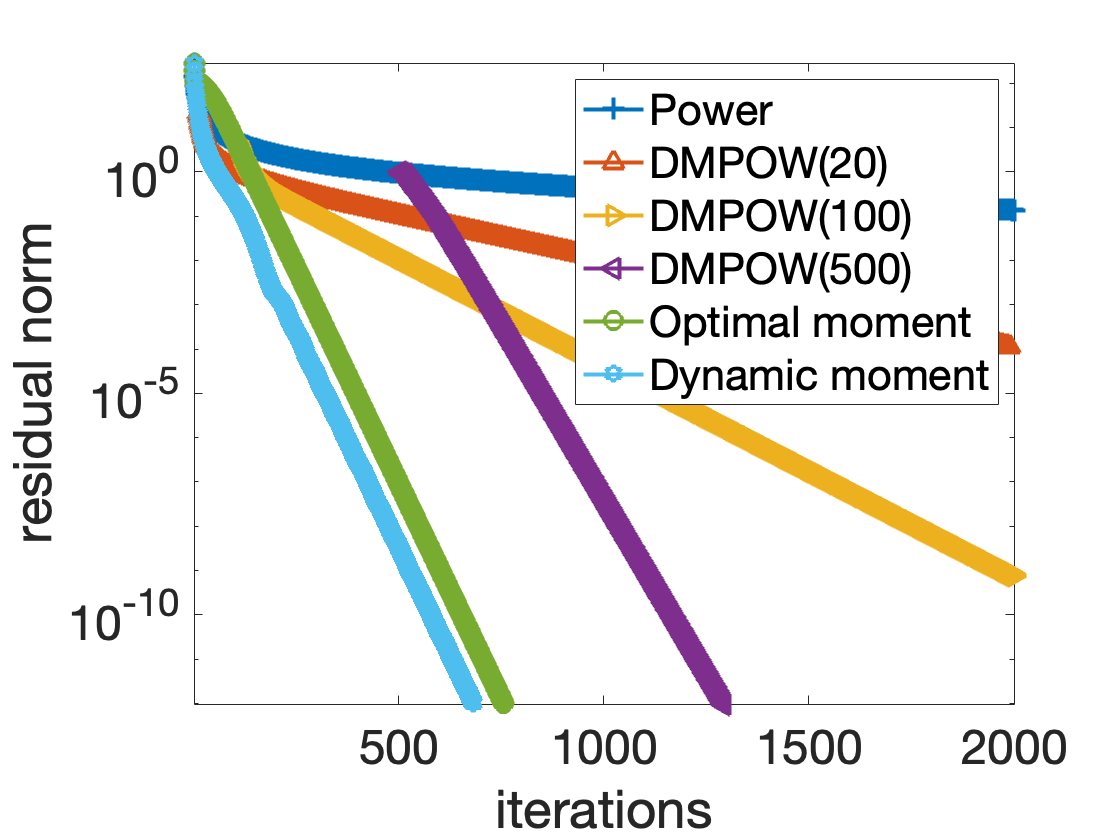}~\hfil~
\includegraphics[trim = 0pt 0pt 0pt 0pt,clip = true, width = 0.32\textwidth]
{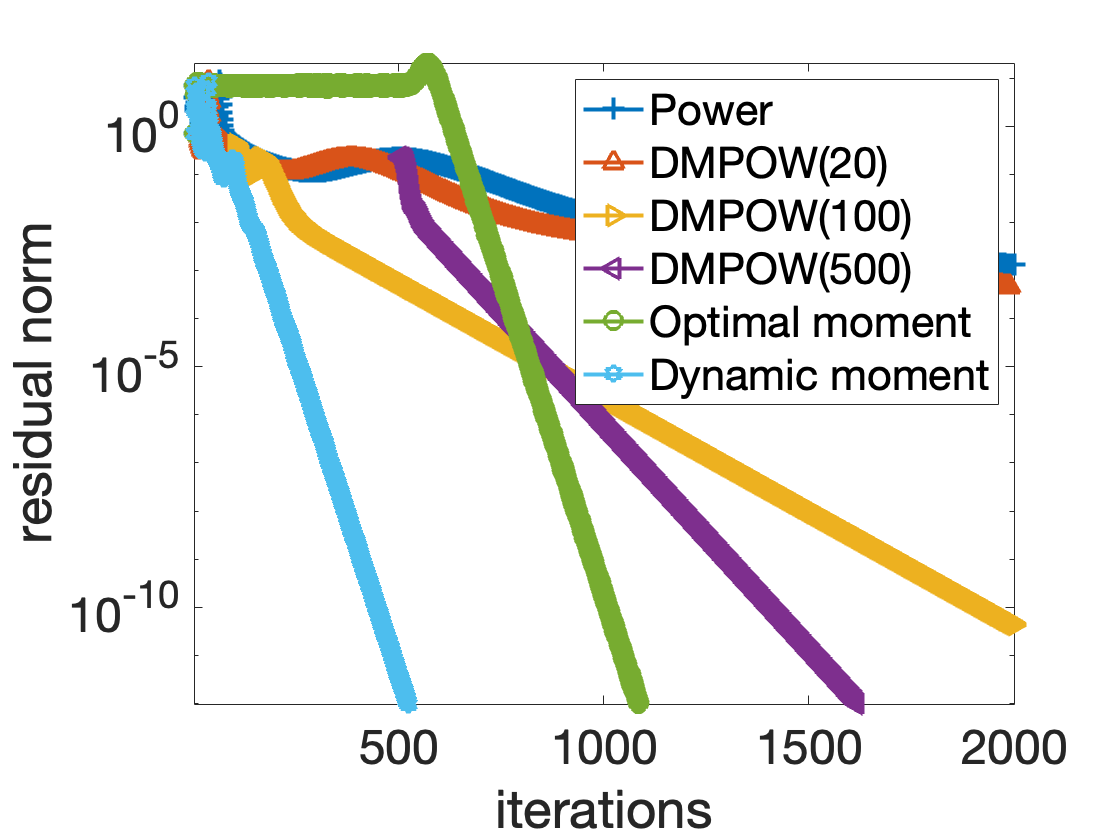}~\hfil~
\includegraphics[trim = 0pt 0pt 0pt 0pt,clip = true, width = 0.32\textwidth]
{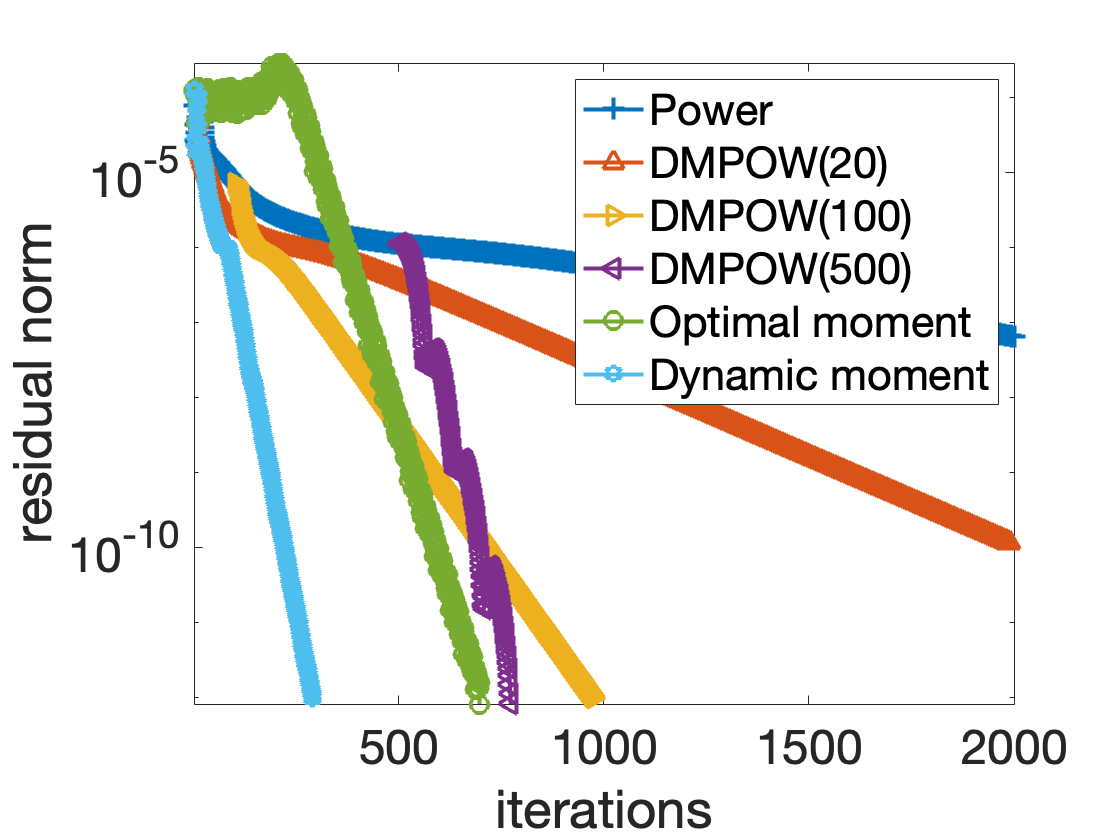}
\caption{Convergence of the residual by iteration count for the three matrices in 
test suite 1, using algorithm \ref{alg:pow}, 
DMPOW with 20, 100 and 500 preliminary iterations,
algorithm \ref{alg:hbpow}, and algorithm \ref{alg:dymo}.
Left: Matrix 1, $\diag(1000:-1:1)$; center: Matrix 2, {\tt Kuu}; 
right: Matrix 3, {\tt Muu}.
\label{fig:test1}
}
\end{figure}

Figure \ref{fig:test1} shows iteration count vs. the residual norm 
using algorithm \ref{alg:pow}, DMPOW with 20, 100 and 500 preliminary iterations,
algorithm \ref{alg:hbpow}, and algorithm \ref{alg:dymo}. Each iteration was
started with the initial $u_0 = \begin{pmatrix}1 & 1 & \cdots & 1 \end{pmatrix}^T$.
The preliminary iterations of DMPOW were started with 
$w_0 = \begin{pmatrix}-1 & 1 & -1 & 1 \cdots \end{pmatrix}^T$, so that $w_0$ is
orthogonal to $u_0$.  In the first two cases, we see the dynamic method 
algorithm \ref{alg:dymo} converges at approximately the same asymptotic rate
as algorithm \ref{alg:hbpow}, though in the second case the latter has an 
extended preasymptotic regime.  The three DMPOW instances work essentially as 
they should for Matrix 1 and Matrix 2, where the approximation of $\lambda_2$
from the deflation method, hence the approximation of $\beta_{opt} = \lambda_2^2/4$
improves as the preliminary iterations are increased. For Matrix 1 DMPOW with 
500 preliminary iterations does appear to achieve the optimal convergence rate.  
In Matrix 3
on the right, only the dynamic method algorithm \ref{alg:dymo} achieves a steady
optimal convergence rate. Algorithm \ref{alg:hbpow} initially stalls then achieves
a good but sub-optimal rate.  DMPOW with 500 preliminary iterations achieves an 
apparently optimal but oscillatory convergence rate, with sub-optimal rates with 100
and 20 preliminary iterations.  The oscillatory behavior of DMPOW suggests that
the approximation to $\lambda_2$ is greater than $\lambda_2$, hence all subdominant
modes are oscillatory, via \eqref{eqn:mucase1}.
%% -- --------------------------------------------------------------------
%% -- --------------------------------------------------------------------
\subsection{Test suite 2}\label{subsec:test2}

The second test suite consists of four matrices. The first three are symmetric 
indefinite and the fourth is SPD with increasing gaps between the smaller eigenvalues. 

\begin{description}
\item[Matrix 4: ] $A = ${\tt ash292} from \cite{DH11}, with $n = 292$. 
This matrix has leading 
eigenvalues $\lambda_1 = 9.1522$ and $\lambda_2 = 8.3769$, with
$r = 0.9153$. It is symmetric indefinite. 
\item[Matrix 5: ] $A = ${\tt bcspwr06} from \cite{DH11}, with $n = 1454$. 
This matrix has leading eigenvalues $\lambda_1 = 5.6195$ and $\lambda_2 = 5.5147$, 
with $r = 0.9814$. It is symmetric indefinite.
\item[Matrix 6: ] $A = \diag(${\tt linspace}$(-99,100,200))$. This matrix has
$n = 200$, $\lambda_1 = 100$, $\lambda_{2,3} = \pm 99$, and $r = 0.99$. It is included to
test the sensitivity to positive and negative leading subdominant eigenvalues.
\item[Matrix 7: ] $A = \diag(10 - ${\tt logspace}$(0,1,200))$.
This matrix has $n = 200$, $\lambda_1 = 9$, $\lambda_2 = 8.9884$, and $r = 0.999$. It is included to test the sensitivity to 
increasing gaps between smaller eigenvalues.
\end{description}

\begin{table}
\centering
\begin{tabular}{|l||c|c||c|c||c|c||c|c|}
\hline
                & \multicolumn{2}{c||}{Matrix 4} &
                  \multicolumn{2}{c||}{Matrix 5} &
                  \multicolumn{2}{c||}{Matrix 6} &
                  \multicolumn{2}{c|}{Matrix 7} \\\hline
method          & min   & max&  min & max 
                & min   & max&  min & max \\\hline
Power     & 247 & 359 & 1088 & 1583& 2000 & 2000& 2000 & 2000 \\\hline
DMPOW(20) & 125 & 746& 197 & 655& 2040 & 2040& 2040 & 2040 \\\hline
DMPOW(100)& 340 & 376& 415 & 1735& 1281 & 2200& 1318 & 2200 \\\hline
DMPOW(500)& 1503 & 1503& 1575 & 1626 & 1767 & 3000& 1970 & 3000 \\\hline
$\beta= \lambda_2^2/4$  & 71 & {\bf 86}& 152 & 179& {\bf 241} & {\bf 288}& 550 & 640 \\\hline
dynamic $\beta$  &  {\bf 66} & 96 & {\bf 133} & {\bf 175}& 255 & 652& {\bf 470} & {\bf 612 }\\\hline 
\end{tabular}
\caption{Minimum and maximum number of matrix-vector multiplies to residual convergence
for 100 runs of the 
power method (algorithm \ref{alg:pow}), 
DMPOW run with 20, 100 and 500 preliminary iterations,
the power iteration with optimal momentum (algorithm \ref{alg:hbpow}) and 
the power iteration with dynamic momentum (algorithm \ref{alg:dymo}), applied
to Matrix 4 - Matrix 7. 
Each run used a randomly generated initial vector.
}
\label{tab:test2}
\end{table}

Results of the experiments with the second set of matrices is shown in table 
\ref{tab:test2}.  We see that the dynamic algorithm \ref{alg:dymo} shows more 
sensitivity to initial vector than does the static algorithm with optimal parameter
\ref{alg:hbpow} in the indefinite cases, and particularly for the highly 
indefinite Matrix 6.
From the results for Matrix 7, we see that increasing the 
spacing between the smaller eigenvalues does not cause increased sensitivity to $u_0$. 
We can also see that algorithms \ref{alg:hbpow} and 
\ref{alg:dymo} significantly outperform the others on all tests in this suite. 

%% -- --------------------------------------------------------------------
%% -- --------------------------------------------------------------------
\subsection{Test suite 3}\label{subsec:test3}
For the third test suite, we generated 100 symmetric matrices with unit
diagonal, and quasi-randomly generated normally distributed off-diagonals
with mean zero and standard deviation one, via
{\tt v = ones(n,1); v1 = randn(n-1,1); A = diag(v,0) + diag(v1,1) + diag(v1,-1);}.
For each matrix, we checked the ratio $r = |\lambda_2/\lambda_1|$. Over the 100
matrices, the values of $r$ ranged from $0.7944$ to $0.9996$, with  mean value $0.9491$ 
and standard deviation $0.0455$. Each run was started with the initial iterate
$u_0 = \begin{pmatrix}1 & 1 & \cdots & 1\end{pmatrix}^T$.

Table \ref{tab:randntest} shows the results.  While algorithm \ref{alg:hbpow} with
optimal fixed $\beta$ has the lowest minimal number of 
iterations over 100 runs, dynamic
algorithm \ref{alg:dymo} has the lowest mean and maximum iteration count.  For these
two methods the iteration count is the same as the reported number of matrix-vector
multiplies.  On the other hand, DMPOW with 20, 100 and 500 preliminary iterations each
had at least one run that did not terminate after 2000 total iterations 
(preliminary included), and all three of the DMPOW methods had a substantially higher
minimum number of matrix-vector multiplies than either the optimal or dynamic methods. 

\begin{table}
\centering
\begin{tabular}{|l||c|c|c|c||c|c|}
\hline
                & \multicolumn{4}{c||}{Matrix-vector multiplies} &
                  \multicolumn{2}{c|}{Terminal residual} \\\hline
method          & mean   &std. dev.&  min& max 
                & min& max \\\hline
Power    & 905.42 & 658.728 & 96 & 2000 & 8.74e-13 & 4.89e-04	\\\hline 
DMPOW(20)&  439.15 & 528.614 & 101 & 2040 & 6.31e-13 & 4.43e-04 \\\hline	 
DMPOW(100)&  498.3 & 295.655 & 302 & 2200 & 2.15e-13 & 7.95e-12 \\\hline	 
DMPOW(500)& 1600.55 & 196.335 & 1502 & 3000 & 1.78e-13 & 4.96e-06 \\\hline
$\beta= \lambda_2^2/4$  & 162.22 & 160.065 & {\bf 40} & 1183 & 6.02e-13 & 9.99e-13 \\\hline
dynamic $\beta$  & {\bf 150.15} & {\bf 133.842} & 63 & \bf{949} & 5.66e-13 & 1.00e-12 \\\hline
\end{tabular}
\caption{The number of matrix-vector multiplies and terminal residual for 100 runs of the 
power method (algorithm \ref{alg:pow}), DMPOW run with 20, 100 and 500 preliminary iterations,
the power iteration with optimal momentum (algorithm \ref{alg:hbpow}) and 
the power iteration with dynamic momentum (algorithm \ref{alg:dymo}). 
Each run used  pseudo-randomly generated tridiagonal matrix
{\tt v = ones(n,1); v1 = randn(n-1,1); A = diag(v,0) + diag(v1,1) + diag(v1,-1);},
and the same initial vector $u_0 = \begin{pmatrix}1 & 1 & \cdots & 1\end{pmatrix}^T$.
}
\label{tab:randntest}
\end{table}

%% -- --------------------------------------------------------------------
%% -- --------------------------------------------------------------------
\subsection{Test suite 4}\label{subsec:test4}
In this fourth suite of tests, we consider three problems 
of varying structure and scale, and which have eigenvalues of varying magnitudes.
The dynamic momentum algorithm \ref{alg:dymo} is tested against 
the power method \ref{alg:pow},
the static momentum
algorithm \ref{alg:hbpow} with optimal parameter $\beta = \lambda_2^2/4$,
and perturbations thereof, 
$\beta_- = 0.99\times\beta_{opt}$, and 
$\beta_+ = \min\{1.01\times\beta_{opt}, (3\lambda_1^2 + \lambda_2^2)/16\}$.
The parameters $\beta_+$ and $\beta_-$ are within 1\%
of $\beta_{opt}$, but do not exceed $\lambda_1^2/4$, which as per section \ref{sec:moment}
would prevent convergence.
\begin{description}
\item[Matrix 8: ] $A = ${\tt Si5H12} from \cite{DH11}, with $n =$ 19,896. 
This matrix has leading eigenvalues 
$\lambda_1 = $58.5609 and $\lambda_2 = $58.4205, with
$r = 0.998$. It is symmetric indefinite. 
\item[Matrix 9: ] $A = ${\tt ss1} from \cite{DH11}, with $n $= 205,282. 
This matrix has leading eigenvalues $\lambda_1 = 1.3735 $ and $\lambda_2 = 1.3733$, 
with $r =0.9998 $. It is nonsymmetric.
\item[Matrix 10: ] $A =${\tt thermomech\_TC}, with $n =$ 102,158. 
This matrix has leading eigenvalues $\lambda_1 = 0.03055$ and $\lambda_2 = 0.03047$, 
with $r =0.9975 $. It is SPD. 
\end{description}

\begin{figure}
\includegraphics[trim = 0pt 0pt 0pt 0pt,clip = true, width = 0.32\textwidth]
{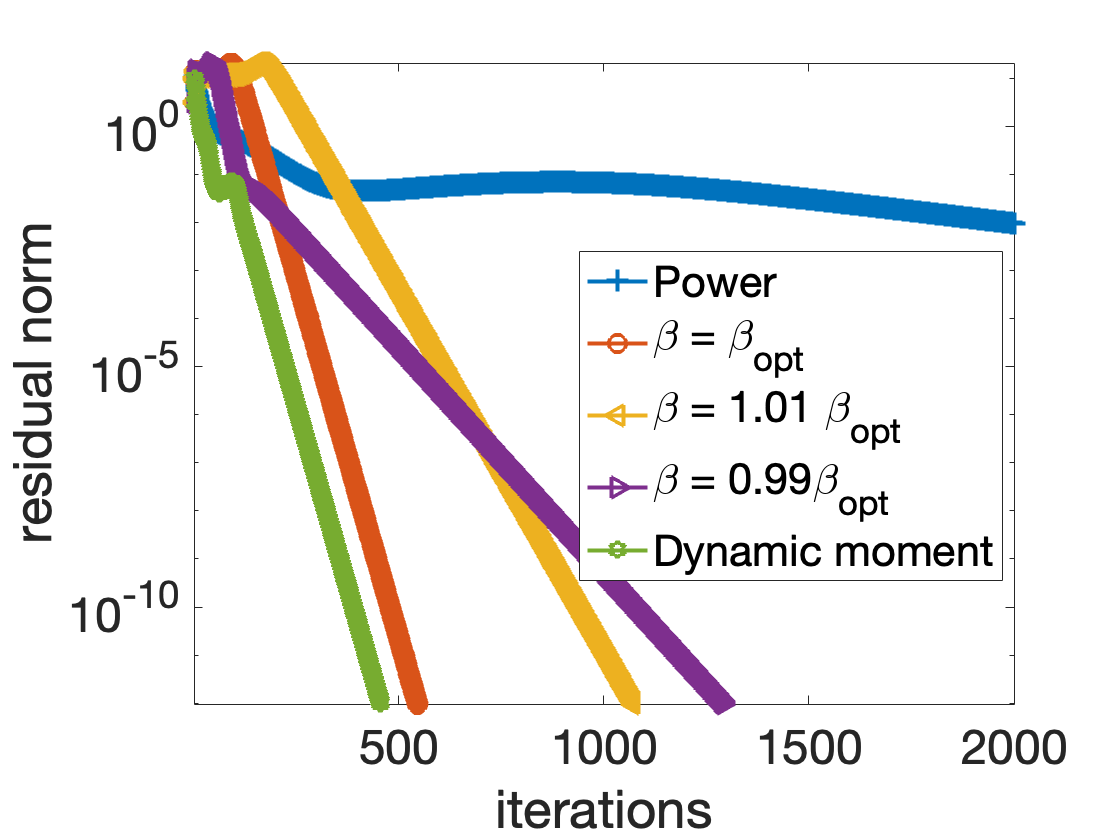}~\hfil~
\includegraphics[trim = 0pt 0pt 0pt 0pt,clip = true, width = 0.32\textwidth]
{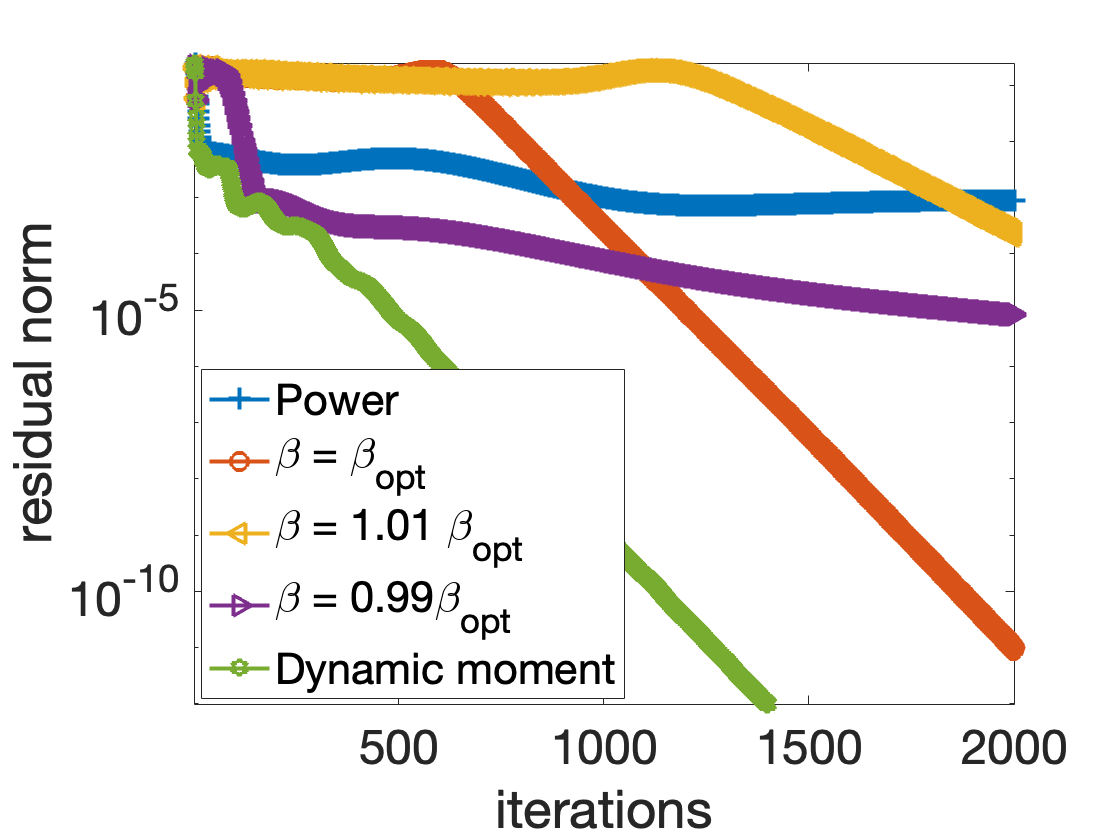}~\hfil~
\includegraphics[trim = 0pt 0pt 0pt 0pt,clip = true, width = 0.32\textwidth]
{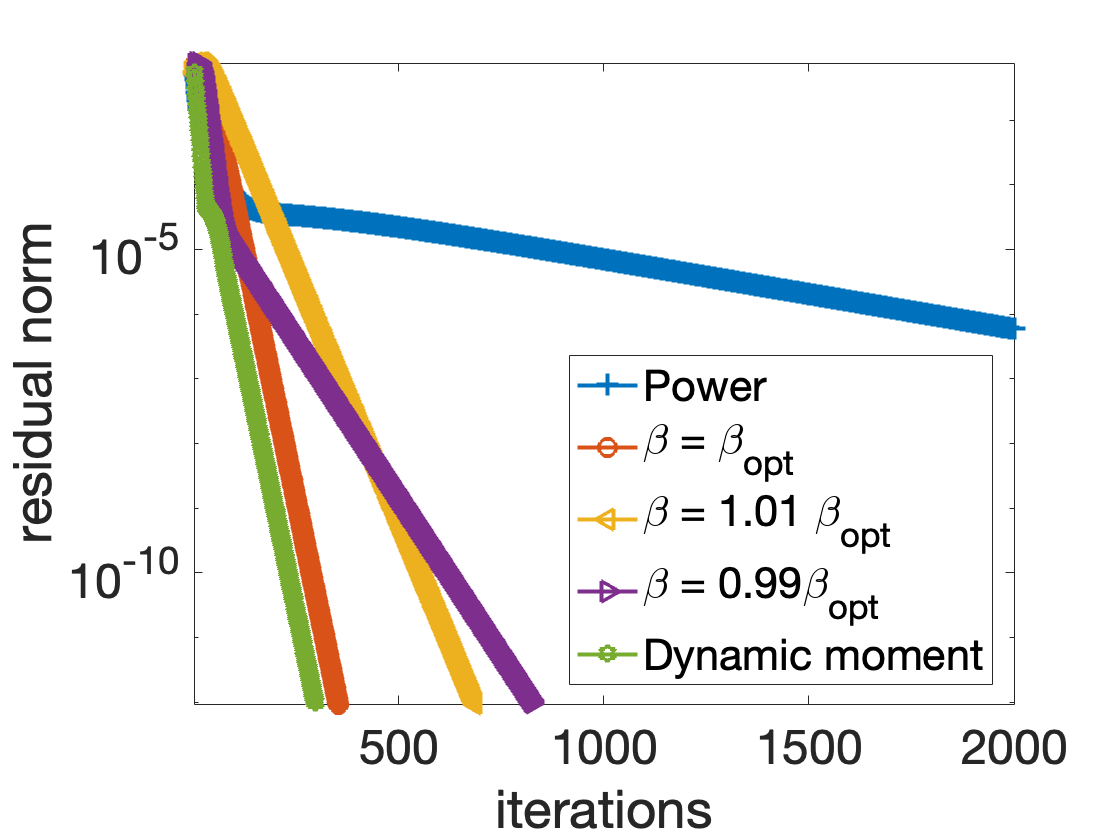}
\caption{Convergence of the residual by iteration count for the three matrices in 
test suite 4, using algorithm \ref{alg:pow}, 
algorithm \ref{alg:hbpow} with $\beta = \beta_{opt} = \lambda_2^2/4$, 
with $\beta = \min\{1.01\times\beta_{opt}, (3\lambda_1^2 + \lambda_2^2)/16\}$,
$\beta = 0.99\times\beta_{opt}$, and 
and algorithm \ref{alg:dymo}.
Left: Matrix 8, {\tt Si5H12}; center: Matrix 9, {\tt ss1}; 
right: Matrix 10, {\tt thermomech\_TC}.
\label{fig:test4}
}
\end{figure}

Convergence of the residual in each case is shown in figure \ref{fig:test4}.
Each of the tests was started from the initial vector 
$u_0 = \begin{pmatrix}1 & 1 & \cdots & 1 \end{pmatrix}^T$.
The results show the dynamic method \ref{alg:dymo} is not sensitive to the scaling
of the eigenvalues which vary in each of the examples.  
The results also show a better rate of convergence with $\beta_+$ than with $\beta_-$, but 
at the cost of a potentially extended preasymptotic regime.  For {\tt ss1} (Matrix 9)
shown in the center plot, the dynamic method shows some initial oscillations but 
does not suffer for the extended asymptotic regime that $\beta_{opt}$ and $\beta_+$
experience. 

%% -- --------------------------------------------------------------------
%% -- --------------------------------------------------------------------
\section{Dynamic Momentum Method for Inverse Iteration}\label{sec:inverse}
%% -- --------------------------------------------------------------------

As an immediate extension of algorithm \ref{alg:dymo}, this section explores the 
application of the dynamic momentum method to accelerate the shifted inverse power 
iteration. 
The shifted inverse iteration is a well-known and powerful technique in the numerical 
solution of eigenvalue problems.
A review of the method including its history, 
theory and implementation can be found in \cite{ipsen97}. 
By appropriately choosing shifting parameters, the inverse iteration with shift can be 
used to identify any targeted eigenpair.
When a good approximation of the targeted eigenvalue is available, the method is
remarkably efficient. 
As the inverse iteration with shift $\sigma$ is equivalent to applying the power iteration on the matrix $(A - \sigma I)^{-1}$ (with the same eigenvectors as those of 
$A$), the analysis carried out in section \ref{subsec:momconv} is directly applicable 
to the algorithm \ref{alg:dymoinv} below.

Each step of the inverse iteration involves solving a linear system. 
For a fixed shift, one can perform a factorization of the shifted matrix before the 
iterative loop to save some computational cost. 
Unlike updating the shift to attain faster convergence, applying a momentum 
acceleration does not require a re-factorization of the matrix.  
As shown below, the momentum accelerated algorithm substantially reduces the number of
iterations to convergence, particularly for sub-optimal shifts. 
Presumably, the use of a sub-optimal shift indicates the user does not have a good
approximation of the target eigenvalue, by which the user is unlikely
to have a good approximation of the second eigenvalue of the shifted system.  Hence
the automatic assignment of the extrapolation parameter $\beta_k$ is essential for
this method to be practical.  Fortunately, as seen below, the proposed method with
dynamic $\beta_k$ is comparable to or outperforms the optimal parameter in each case
tested.
Numerical experiments  below illustrate the improved efficiency, particularly with 
the dynamic strategy.

In our implementation of DMPOW in this section we terminated the preliminary iterations
in the deflation stage when the target (second) eigenvalue achieved a given 
relative tolerance, that is in the notation of \cite{RJRH22}, 
$|(\mu_j - \mu_{j-1})/\mu_j| < 10^{-n}$. We show results
using $n = 1,2,4$.  We implemented DMPOW as a shifted inverse iteration as follows,
referring to the implementation given in \cite[algorithm 1]{RJRH22}:
We replaced the multiplication by $A$ in line 2 with a multiplication by 
$(A - \sigma I)^{-1}$, implemented as a solve of the LU-factored system, and the
multiplication by $A-P$ in line 6 with a multiplication by $(A-\sigma I)^{-1}$, again
implemented as a solve of the factored system, and a multiplication by $P$.  We
replaced the multiplication by $A$ in line 8 to compute the Rayleigh quotient 
corresponding to the second eigenvalue as a multiplication by $A - \sigma I$.
We did this rather than multiplying by $(A - \sigma I )^{-1}$ to reduce the number of 
system solves per iteration from three to two, with little to no effect on the total
iteration count.

Just as algorithm \ref{alg:dymo} requires two preliminary power iterations, the
dynamic momentum strategy for the inverse iteration requires two preliminary 
inverse iterations.  
For tests in this section, we ran iterations to a residual tolerance of 
$10^{-15}$ or a maximum of 2000 iteration.
In this section we also numerically verify the stability of the extrapolation
parameter $\beta_k$ as shown in lemma \ref{lem:rstab} and remark \ref{rem:betastab}.
%% -- --------------------------------------------------------------------
\begin{algorithm}{Inverse power iteration}
\label{alg:ipow}
\begin{algorithmic}
\State{Choose $v_0$ and shift $\sigma$, set $h_0 = \nr{v_0}$ and $x_0 = h_0^{-1}v_0$}
\State{Compute $(A- \sigma I) = LU$} \Comment{Compute LU factors of $A-\sigma I$}
%\State{Set $v_1 = A x_0$}
\State{Solve $L y = x_0$ and  $U v_{1} = y$} 
\For{$k \ge 0$}
\State{Set $h_{k+1} = \nr{v_{k+1}}$ and  $x_{k+1} = h_{k+1}^{-1} v_{k+1}$}
\State{Solve $L y = x_{k+1}$ and  $U v_{k+2} = y$} 
\State{Set $\nu_{k+1} = (v_{k+2}, x_{k+1})$ and 
       $d_{k+1} = \|v_{k+2} - \nu_{k+1} x_{k+1}\|$}
\State{{STOP} if $\nr{d_{k+1}} <$ \tt{tol}  }
\EndFor
\end{algorithmic}
\end{algorithm}
%% -- --------------------------------------------------------------------

The dynamically accelerated version of algorithm \ref{alg:ipow} follows.
%% -- --------------------------------------------------------------------
\begin{algorithm}{Dynamic momentum for inverse iteration}
\label{alg:dymoinv}
\begin{algorithmic} 
\State {Do two iterations of algorithm \ref{alg:ipow}} \Comment{$k = 0,1$}
\State{Set $r_{2} = \min\{d_2/d_{1},1\}$}
\For{$k \ge 2$} \Comment{$k \ge 2$}
\State{Set $\beta_k = \nu_{k}^2 r_{k}^2/4$}
\State{Set $u_{k+1} = v_{k+1} - (\beta_k/h_k) x_{k-1}$}
\State{Set $h_{k+1} = \nr{u_{k+1}}$ and $x_{k+1} = h_{k+1}^{-1} u_{k+1}$}
\State{Solve $L y = x_{k+1}$ and  $U v_{k+2} = y$} 
\State{Set $\nu_{k+1} = (v_{k+2}, x_{k+1})$ and
           $d_{k+1} = \nr{v_{k+2} - \nu_{k+1} x_{k+1}}$}
\State{Update $\rho_k = \min\{d_{k+1}/d_k,1\}$  and
              $r_{k+1} = 2\rho_k/(1 + \rho_k^2)$}
\State{{STOP} if $\nr{d_{k+1}} <$ \tt{tol}  }
\EndFor
\end{algorithmic}
\end{algorithm}
%% -- --------------------------------------------------------------------
%% -- --------------------------------------------------------------------

In table \ref{tab:inv} we show results for accelerating the inverse
iteration used to recover the largest eigenvalue of the matrix $A = \diag(1000:-1:1)$.
We test shifts $\sigma = \{999.75,\ldots,1064\}$,
chosen with increasing distance from the target eigenvalue $\lambda_1 = 1000$
to see how much a suboptimal shift can be made up for with the extrapolation. 
We see the dynamic momentum method and the ``optimal'' fixed momentum parameter 
give the best performance, with the dynamic method converging in fewer iterations
as the shift increases away from the target eigenvalue.  
The performance of the DMPOW iteration is intermediate
between the base inverse iteration and the dynamical momentum method.
Compared with the base algorithm \ref{alg:ipow},
algorithm \ref{alg:dymoinv} with dynamically chosen $\beta_k$ 
not only reduces the number of iterations for each given shift, it also achieves a better 
iteration count with shifts more than twice as far away from the target eigenvalue.  
This shows algorithm
\ref{alg:dymoinv} reduces the sensitivity to the shift in the standard inverse iteration.
\begin{table}
\centering
\begin{tabular}{|l||c||c||c||c|c|c|}%c|c|c|c|}
\hline
~$\sigma$ 
& $\beta = 0$ & \bf{dyn $\beta_k$} & $\beta_{opt}$ & 
$\beta_{DM}(10^{-1})$ & $\beta_{DM}(10^{-2})$ & $\beta_{DM}(10^{-4})$ \\\hline
\ 999.75 & 33 & {\bf 21}  & 23 (4.44e-1) & 25 (4.60e-1) & 29 (4.45e-1) & 32 (4.44e-1) \\\hline
1000.25  & 23 & {\bf 17}  & {\bf 18} (1.60e-1) & 21 (1.59e-1)& 22 (1.60e-1) & 25 (1.60e-1) \\\hline
1000.5   & 32 & {\bf 23}  & {\bf 22} (1.11e-1) & 26 (1.11e-1)& 26 (1.11e-1) & 31 (1.11e-1) \\\hline
1001     & 49 & 33  & {\bf 29} (6.25e-2)& 32 (6.39e-2)& 32 (6.39e-2) & 40 (6.25e-2) \\\hline
1004     & 142 & 55  & {\bf 52} (1.00e-2)& 56 (1.02e-2) & 58 (1.03e-2) & 80 (1.00e-2) \\\hline
1016     & 478 & {\bf 88}  & 95 (8.65e-4)& 215 (7.43e-4) & 134 (8.52e-4) & 120 (8.73e-4) \\\hline
1064     & 1691 & {\bf 163} & 175 (5.92e-5) & 841 (4.55e-5) & 525 (5.47e-5)& 239 (5.93e-5) \\\hline
\end{tabular}
\caption{
Number of system solves using the shifted inverse iteration to find the largest 
eigenvalue  with fixed and dynamic momentum for the matrix
$A = \diag(n:-1:1)$, with $n=1000$.  The optimal fixed extrapolation parameters
$\beta_{opt} = 1/(4(\lambda_2- \sigma)^2)$ are shown after the number of solves
in the $\beta_{opt}$ column, and the dynamically chosen parameter $\beta_k$ is set
as in algorithm \ref{alg:dymoinv}. 
The last three columns, $\beta_{DM}(10^{-n}), n = 1,2,4$ contain the number of system 
solves for DMPOW, where the 
preliminary deflation stage is terminated after the target eigenvalue reaches the
relative tolerance of $10^{-n}$. The parameter is shown after the number of solves.
Each iteration is started from 
$v_0 = \begin{pmatrix} 1 & 1 &\cdots& 1\end{pmatrix}^T$, and
run to a residual tolerance of $10^{-15}$.
\label{tab:inv}
}
\end{table}
\begin{figure}
\includegraphics[trim = 0pt 0pt 0pt 0pt,clip = true, width = 0.32\textwidth]
{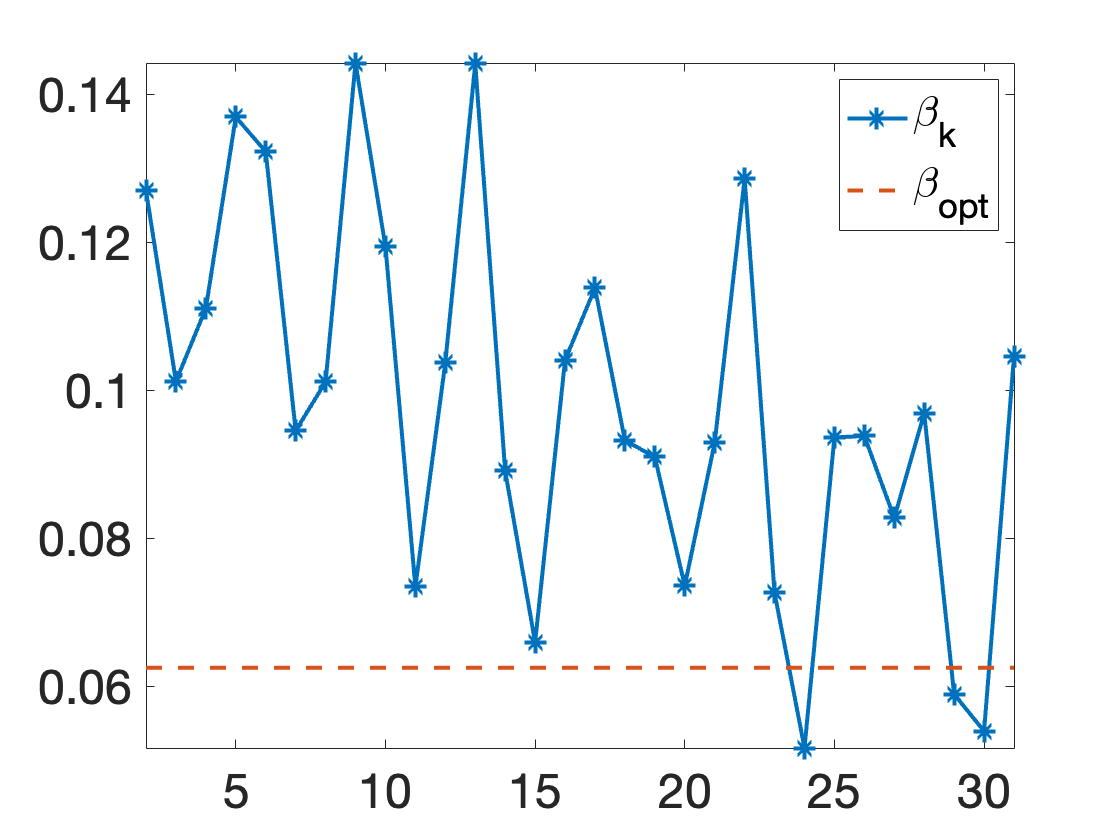}~\hfil~
\includegraphics[trim = 0pt 0pt 0pt 0pt,clip = true, width = 0.32\textwidth]
{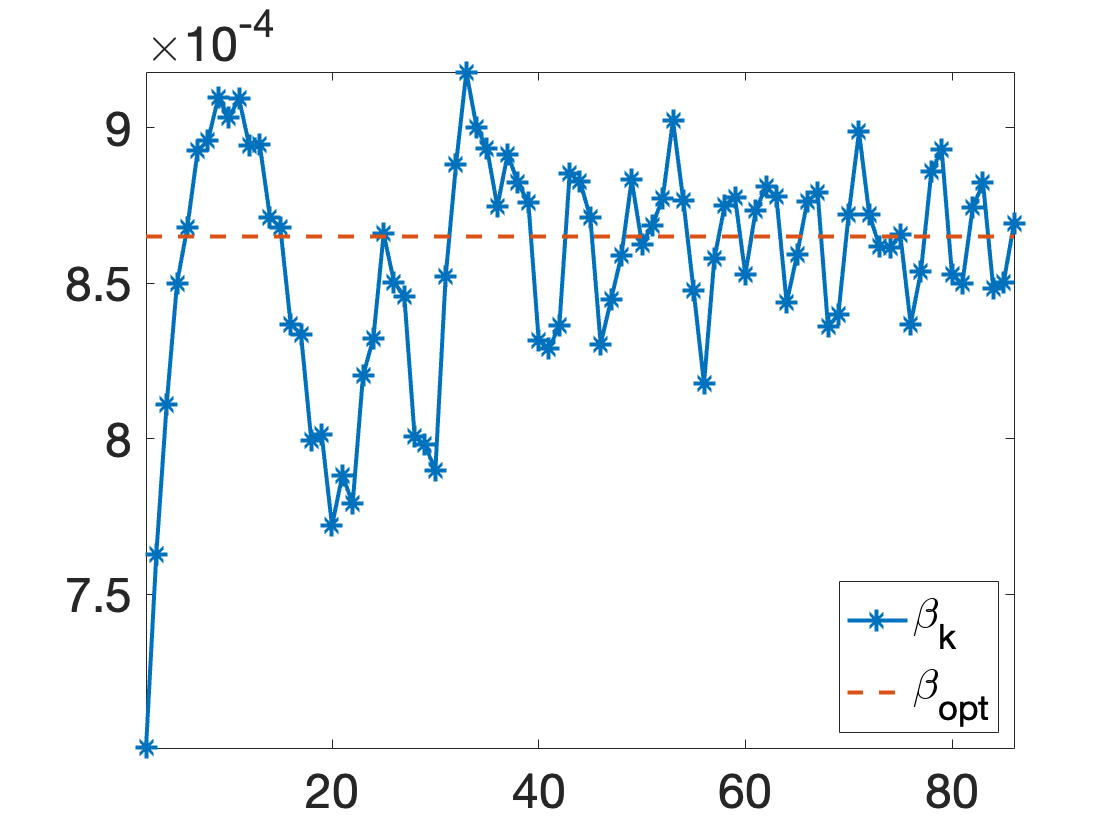}~\hfil~
\includegraphics[trim = 0pt 0pt 0pt 0pt,clip = true, width = 0.32\textwidth]
{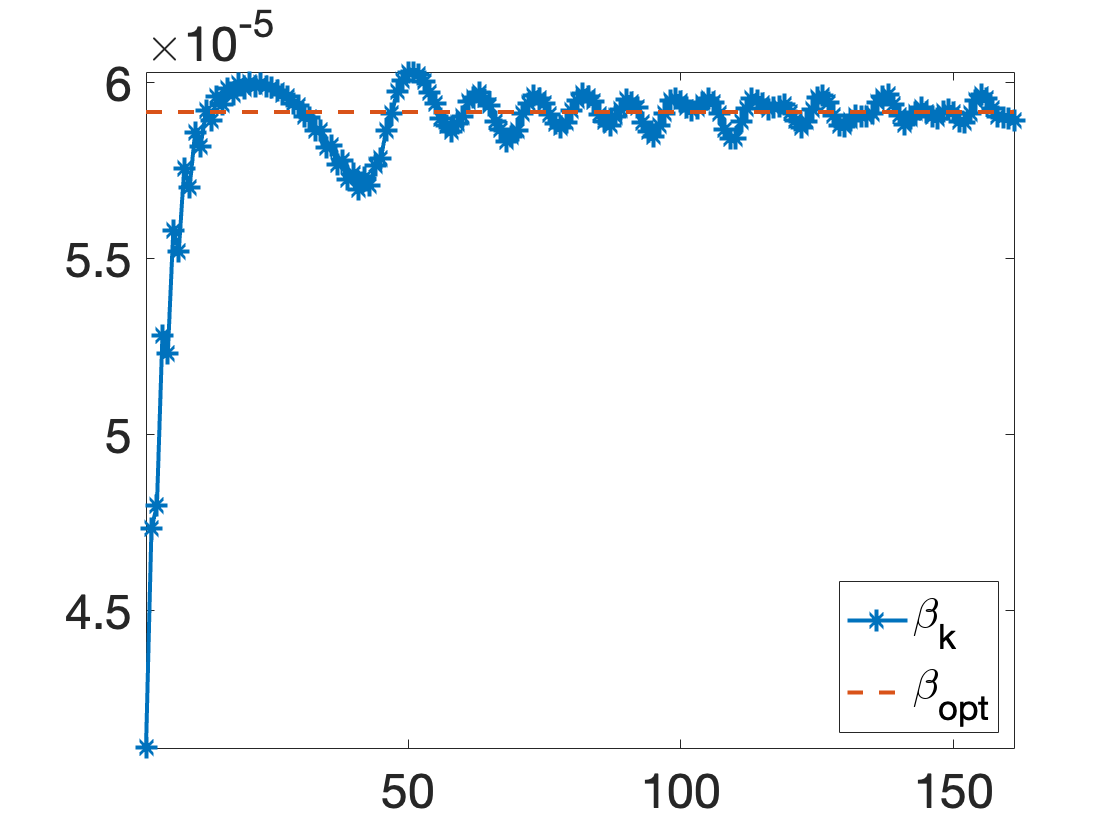}
\caption{Behavior of $\beta_k$ with respect to $\beta_{opt}$
for representative examples from table \ref{tab:inv}, 
illustrating that $\beta_k$ stabilizes closer
to $\beta_{opt}$ in agreement with lemma \ref{lem:rstab} as $r \goto 1$. 
Left: $\sigma = 1001$, for which $r = 0.5$.
Center: $\sigma = 1016$, for which $r \approx 0.94$.  Right: $\sigma = 1064$ for which 
$r \approx 0.98$. 
\label{fig:invpmbetas}
}
\end{figure}

Figure \ref{fig:invpmbetas} 
shows the extrapolation parameters $\beta_k$ for three different
shifts $\sigma$ as shown in table \ref{tab:inv}.  The left plot shows $\{\beta_k\}$ for 
$\sigma = 1001$.  Denoting the eigenvalues of $(A - \sigma I )^{-1}$ as 
$\{\tilde \lambda_i\}, i = 1, \ldots, n$, we have
$\tilde \lambda_1 = -1$, and $\tilde \lambda_2= -1/2$ so that $r = 0.5$. 
In this case $\beta_k$ oscillates above and below $\beta_{opt}$.  
For $\beta_k < \beta_{opt}$, $\tilde \lambda_2$ is non-oscillatory by \eqref{eqn:mucase1}, 
but each of the eigenvalues with $\tilde \lambda ^2/4 < \beta$ is oscillatory, 
and each decays at a 
slightly faster rate than $\tilde \lambda_2$ by \eqref{eqn:mucase3}.  
As per lemma \ref{lem:rstab}, the approximation of $r$ by the detected convergence 
rate $\rho_k$ is stable, but the oscillations in $r_{k+1}$ are not necessarily 
damped with respect to the detected $\rho_k$.
For $\beta_k > \beta_{opt}$, all subdominant modes are oscillatory and decay at the same
rate by \eqref{eqn:mucase3}, and the stability of $r_{k+1}$ with respect to $\rho_k$ 
still holds.
For the center plot in figure~\ref{fig:invpmbetas}, $\sigma = 1016$ so that 
$r \approx 0.94$; and in the right plot $\sigma = 1064$ so that $r \approx 0.98$. 
In both of these cases, lemma \ref{lem:rstab}
shows that $r_{k+1}$ is stable with respect to $\rho_k$, and the difference between $r_k$ and
$r$ is damped in comparison to the difference between the detected $\rho_k$ and $\rho$; 
and moreso in the plot on the right.  
Notably for the center figure with $r \approx 0.94$, $\beta_k$ converges to $\beta_{opt}$
to within $10^{-4}$, and for the right plot with $r \approx 0.98$, $\beta_k$ 
converges to $\beta_{opt}$ to within $10^{-5}$.
This demonstrates how $r_k$ approaches $r$ as $r$ approaches one, as described in lemma 
\ref{lem:rstab}.    

\begin{table}
\centering
\begin{tabular}{|l||c||c||c||c|c|c|c|c|c|c|}
\hline
~$\sigma$ 
& $\beta = 0$ & \bf{dyn $\beta_k$} & $\beta_{opt}$ & 
$\beta_{DM}(10^{-1})$ & $\beta_{DM}(10^{-2})$ & $\beta_{DM}(10^{-4})$ \\\hline
\ \ 1.25 & 33 & 21 & 23 (4.44e-1) & 25 (4.60e-1) & 29 (4.45e-1) & 32 (4.44e-1)\\\hline 
\ \ 0.75 & 23 & 17 & 17 (1.60e-1) & 21 (1.59e-1) & 22 (1.60e-1) & 25 (1.60e-1)\\\hline 
\ \ 0    & 49 & 33 & 29 (6.25e-2) & 32 (6.39e-2)& 32 (6.39e-2) & 40 (6.25e-2) \\\hline 
\ \ -1   & 81 & 46 & 39 (2.78e-2) & 41 (2.89e-2)& 42 (2.86e-2) & 57 (2.78e-2)\\\hline 
\ \ -4   & 171 & 58 & 57 (6.94e-3) & 59 (6.97e-3)& 64 (7.10e-3) & 89 (6.95e-3) \\\hline
\ \ -8   & 286 & 70 & 74 (2.50e-3) & 119 (2.32e-3) & 78 (2.52e-3) & 120 (2.50e-3) \\\hline
\ \ -16  & 505 & 91 & 97 (7.72e-4) & 229 (6.58e-4) & 143 (7.57e-4) & 124 (7.79e-4)\\\hline
\ \ -32  & 922 & 123&130 (2.16e-4) & 444 (1.73e-4) & 288 (2.03e-4)& 177 (2.17e-4) \\\hline 
\end{tabular}
\caption{
Number of system solves using the shifted inverse iteration to find the smallest 
eigenvalue  with fixed and dynamic momentum for the matrix
$A = \diag(n:-1:1)$, with $n=1000$.  The optimal fixed extrapolation parameters
$\beta_{opt} = 1/(4(\lambda_{n-1}- \sigma)^2)$ are shown after the number of solves
in the $\beta_{opt}$ column, and the dynamically chosen parameter $\beta_k$ is set
as in algorithm \ref{alg:dymoinv}. 
The last three columns, $\beta_{DM}(10^{-n}), n = 1,2,4$ contain the number of system 
solves for DMPOW, where the 
preliminary deflation stage is terminated after the target eigenvalue reaches the
relative tolerance of $10^{-n}$. The parameter is shown after the number of solves.
Each iteration is started from 
$v_0 = \begin{pmatrix} 1 & 1 &\cdots& 1\end{pmatrix}^T$, and
run to a residual tolerance of $10^{-15}$.
\label{tab:invsmall}
}
\end{table}

\begin{figure}
\includegraphics[trim = 0pt 0pt 0pt 0pt,clip = true, width = 0.32\textwidth]
{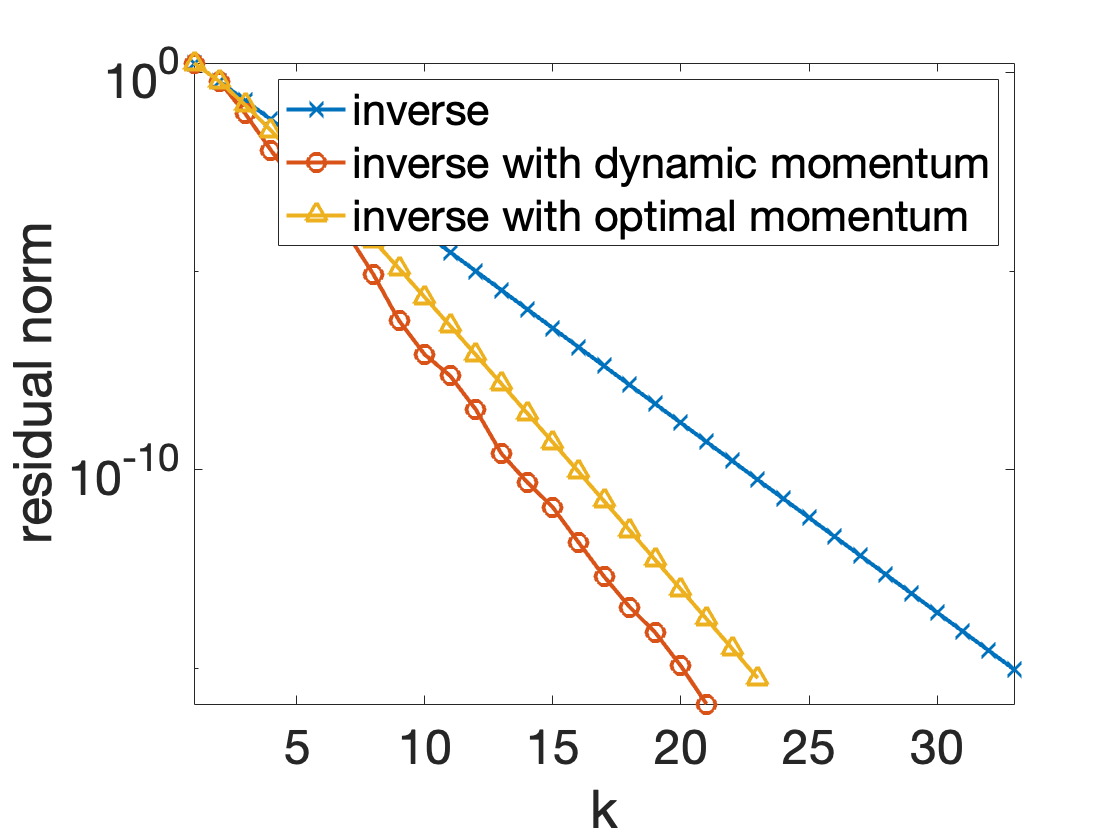}~\hfil~
\includegraphics[trim = 0pt 0pt 0pt 0pt,clip = true, width = 0.32\textwidth]
{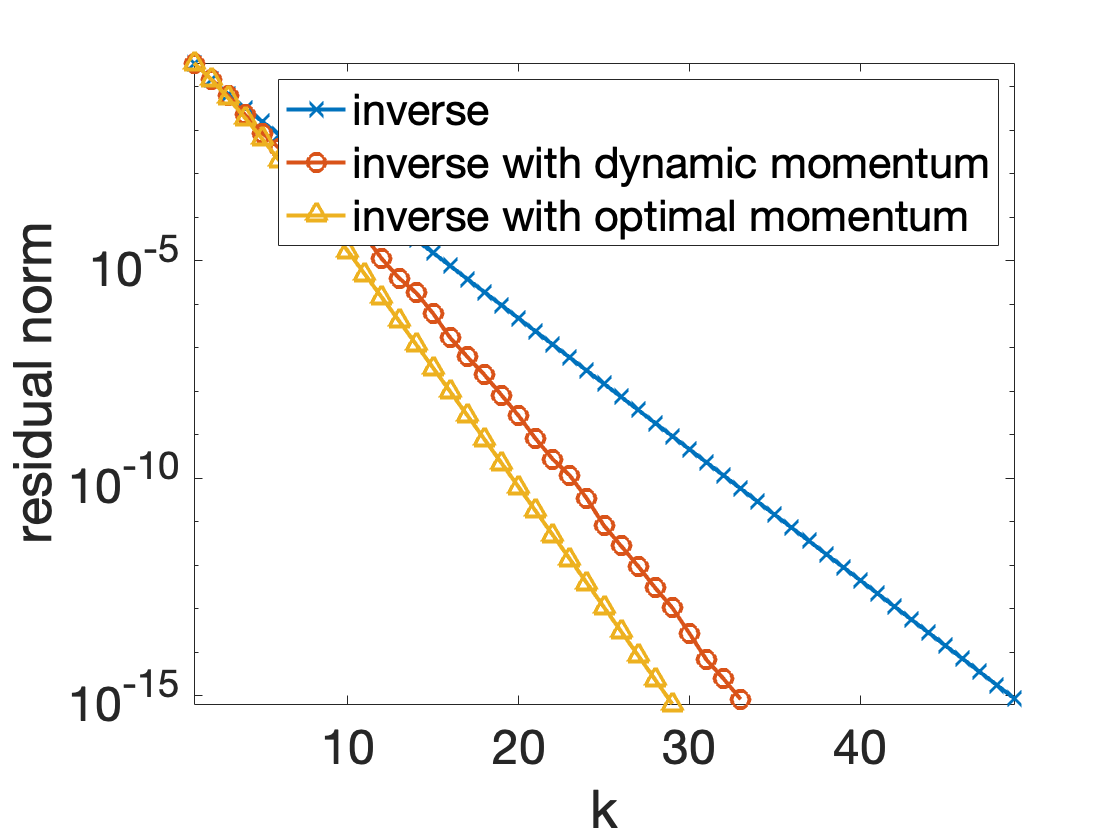}~\hfil~
\includegraphics[trim = 0pt 0pt 0pt 0pt,clip = true, width = 0.32\textwidth]
{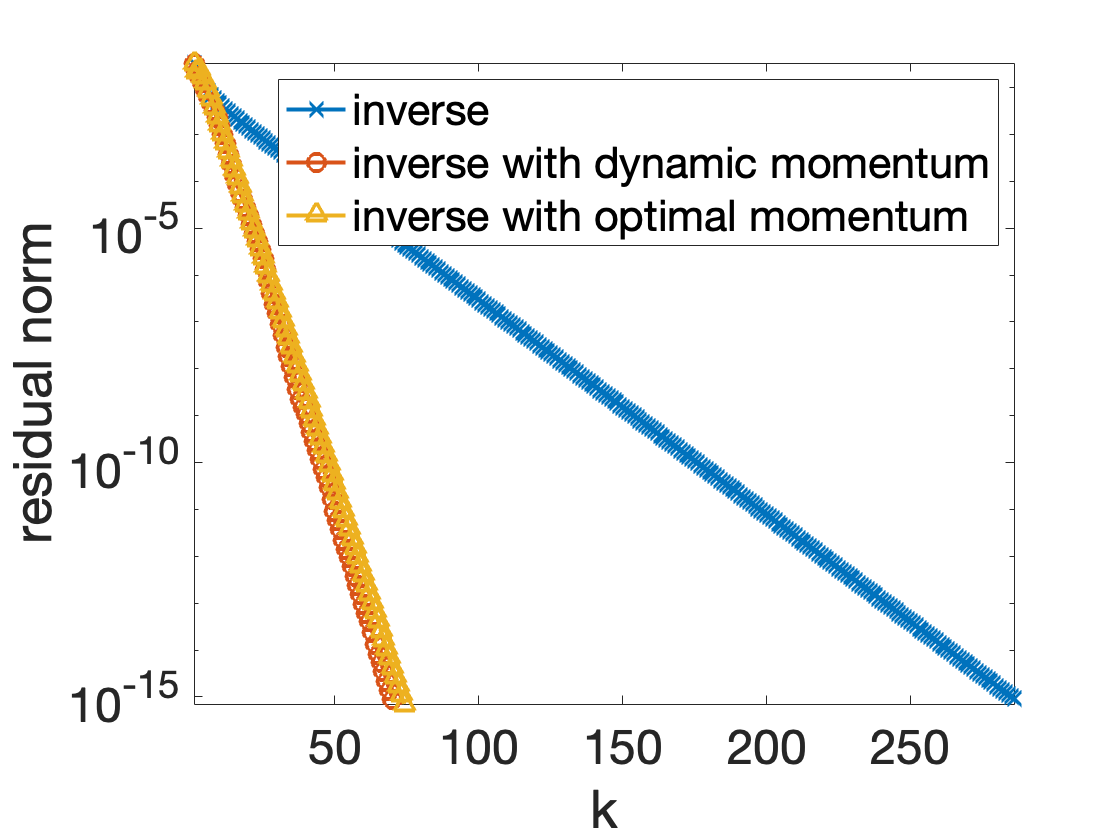}
\caption{Residual convergence for representative examples from table \ref{tab:invsmall}, 
illustrating the improvement in convergence for the optimal and dynamic momentum 
acceleration for a variety of shifts.
Left: $\sigma = 1.25$, for which $r = 1/3$.
Center: $\sigma = 0$, for which $r = 0.5$.  Right: $\sigma = -8$ for which 
$r = 0.9$. 
\label{fig:invpmconv}
}
\end{figure}

In table \ref{tab:invsmall} we show the results of a similar experiment to recover
the smallest eigenvalue of the matrix $A = \diag(1000:-1:1)$.
We test shifts 
$\sigma = \{1.25, 0.75, 0.5, 0, -1, -4, -8, -16, -32\}$, 
a range of shifts with increasing distance from the target eigenvalue $\lambda_n = 1$.
Our results are similar to the largest eigenvalue case of table \ref{tab:inv}.  
We see the ratio between the number of iterations in the dynamic method and the 
inverse iteration decreases as the the shift increases. For instance, as the shift
ranges between 0 and -32, the corresponding ratio between the number of dynamic
momentum iterations and inverse iterations without momentum decreases monotonically
from 0.67 to 0.13.
In each case tested, the dynamic method is either comparable to or better than
the momentum method with optimal shift, and outperforms all of the DMPOW iterations. 
Figure \ref{fig:invpmconv} shows convergence plots for $\sigma = 1.25, \sigma = 0$ 
and $\sigma = -8$, providing a visualization of 
the improved convergence rates from the dynamic algorithm \ref{alg:dymoinv}.
These examples illustrate the gain in convergence from
this practical and low-cost acceleration method.

%% -- --------------------------------------------------------------------
%% -- --------------------------------------------------------------------
\section{Conclusion}\label{sec:conc}
%% -- --------------------------------------------------------------------
%% -- --------------------------------------------------------------------
In this paper we introduced and analyzed a one step extrapolation method to accelerate 
convergence of the power iteration for real, diagonalizable matrices, 
and proved convergence to the dominant eigenpair with acceleration in the symmetric case.
The method is based on the momentum method for the power iteration introduced
in \cite{DSHMRX19}, 
and requires a single matrix-vector multiply per iteration. 
Unlike the method of \cite{DSHMRX19} and other recent variants 
such as \cite{RJRH22}, the presently introduced technique gives a dynamic update of
the key extrapolation parameter at each iteration, and does not require any 
a priori knowledge of the spectrum.

We first reviewed some results on the analysis of a static method of the 
form introduced in \cite{DSHMRX19}, by considering the power iteration 
applied to an augmented matrix. Our analysis goes beyond that shown in the
original paper, revealing that the augmented matrix is defective for the 
optimal parameter choice, which explains why slower convergence is expected
in the preasymptotic regime.  We then analyzed our dynamic method showing 
both stability of the dynamic extrapolation parameter and convergence of the method.

In the last two sections we numerically demonstrated the efficiency of the introduced 
dynamic algorithm \ref{alg:dymo} as applied to power and inverse iterations.
We demonstrated that algorithm \ref{alg:dymo} 
often outperforms the original static method with the optimal parameter choice
as given in algorithm \ref{alg:hbpow}. We further showed algorithm \ref{alg:dymo}
performs favorably in comparison to the method of \cite{RJRH22}, which generally 
accelerates the power iteration but does not exceed the performance of 
algorithm \ref{alg:hbpow}. Finally, we showed that the introduced dynamic method 
is a useful tool to
accelerate inverse power iterations, and can be used to converge in as few iterations as 
having
a shift twice as close to the target eigenvalue, and without significant additional
computational complexity.
Future work will include the development of an analogous method applied to 
(preconditioned) Krylov subspace projection methods as in 
\cite{GoYe02,knyazev01,PoSc21,QuYe10} to efficiently recover multiple eigenpairs.
%% -- --------------------------------------------------------------------
%% -- --------------------------------------------------------------------
\section{Acknowledgements}\label{sec:ack}
CA and SP are supported in part by NSF DMS-2045059 (CAREER).
This material is based upon work supported by the NSF under DMS-1929284 while 
SP was in residence at the Institute for Computational and Experimental Research in
Mathematics in Providence, RI, during the Numerical PDEs: Analysis, Algorithms and 
Data Challenges Program.  SP would like to thank Prof. Nilima Nigam for many 
interesting discussions that led to the formulation of this work.
The authors would also like to thank the anonymous referees, whose thoughtful 
suggestions improved this manuscript.
%% -- --------------------------------------------------------------------
%% -- --------------------------------------------------------------------

%% -- --------------------------------------------------------------------
%% -- --------------------------------------------------------------------

\bibliographystyle{abbrv}
\bibliography{accrefs}
%% -- --------------------------------------------------------------------
%% -- --------------------------------------------------------------------

\end{document}